\documentclass[12pt]{article}
\usepackage[margin=2.3cm,a4paper]{geometry}
\usepackage{amsmath}
\usepackage{amssymb}
\usepackage{amsthm}
\usepackage{graphicx,epstopdf}
\usepackage[algo2e,vlined,longend,linesnumbered,ruled]{algorithm2e}
\SetCommentSty{itshape}
\SetKwComment{Comment}{}{}
\usepackage{tikz}
\usetikzlibrary{positioning,shapes,calc,backgrounds,fit,plotmarks}
\usepackage{pgfplots}
\usepackage{subfig}
\usepackage{textcomp} 
\usepackage{xspace}
\usepackage{afterpage}
\usepackage{hyperref}

\pgfplotsset{compat=newest,
every axis/.append style={axis x line=bottom,
                          axis y line=left,
                          scale only axis,
                          y label style={at={(0.1,1.0)},anchor=south west,rotate=-90,color=black},
                          },
}

\theoremstyle{definition}

\newcommand\R{\ensuremath{\mathbb{R}}}
\newcommand{\cA}{\ensuremath{\mathcal{A}}}
\newcommand{\cB}{\ensuremath{\mathcal{B}}}
\newcommand{\cC}{\ensuremath{\mathcal{C}}}

\newcommand{\cL}{\ensuremath{\mathcal{L}}}
\newcommand{\vvec}{\mathrm{vec}}
\newcommand{\trace}[1]{\ensuremath{\operatorname{Tr}\!\left(#1\right)}}
\newcommand{\micron}{{\textmu}m\xspace}
\newcommand{\dee}{\mathrm{d}}
\newcommand{\trans}{\intercal} 
\newcommand{\rank}{r}
\newcommand{\ninj}{{n_\text{inj}}}
\newcommand{\nx}{{n_\text{X}}}
\newcommand{\ny}{{n_\text{Y}}}
 \newcounter{mymac@matlab}
\newcommand{\matlab}{MATLAB%
   \ifnum\value{mymac@matlab}<1%
   \textsuperscript{\textregistered}%
   \setcounter{mymac@matlab}{1}%
   \fi%
  }
\newcommand{\intel}{Intel\textsuperscript{\textregistered}}
\newcommand{\nvidia}{{\fontfamily{jkpss}\fontshape{sc}\selectfont Nvidia}\textsuperscript{\textregistered}}

\usepackage{tikz} \usetikzlibrary{positioning,shapes,calc,backgrounds,fit,plotmarks}
\usepackage{pgfplots}
\pgfplotsset{compat=newest,
every axis/.append style={axis x line=bottom,
                          scale only axis,
                          },
}

\date{\today}

\begin{document}
\title{Greedy low-rank algorithm for spatial connectome regression}
\author{%
  Patrick K{\"u}rschner\footnote{KU Leuven, Electrical Engineering (ESAT/STADIUS), Campus Kulak Kortrijk 
	(\texttt{patrick.kurschner@kuleuven.be}).}
  \and Sergey Dolgov\footnote{University of Bath, Mathematical Sciences (\texttt{s.dolgov@bath.ac.uk})}
  \and Kameron Decker Harris\footnote{University of Washington, Computer Science \& Engineering (\texttt{kamdh@uw.edu}).}
  \and Peter Benner\footnote{Max Planck Institute for Dynamics of Complex Technical Systems, Magdeburg (\texttt{benner@mpi-magdeburg.mpg.de}).}
}

\maketitle

\begin{abstract}
Recovering brain connectivity from tract tracing data is an
important computational problem in the neurosciences.
Mesoscopic connectome reconstruction
was previously formulated as a structured matrix regression problem \cite{Harris16},
but existing techniques do not scale to the whole-brain setting.
The corresponding matrix equation is challenging to solve due to large scale,
ill-conditioning, and a general form that lacks a convergent splitting.
We propose a greedy low-rank algorithm for
connectome reconstruction problem in very high dimensions.
The algorithm approximates the solution by a sequence
of rank-one updates which exploit the sparse and positive definite problem structure.
This algorithm was described previously \cite{KreS15}
but never implemented for this connectome problem,
leading to a number of challenges.
We have had to design judicious stopping criteria
and employ efficient solvers for the three main sub-problems of the algorithm,
including an efficient GPU implementation 
that alleviates the main bottleneck for large datasets.
The performance of the method is evaluated on three examples:
an artificial ``toy'' dataset and
two whole-cortex instances using
data from the Allen Mouse Brain Connectivity Atlas.
We find that the method is significantly faster than previous
methods and that moderate ranks offer good approximation.
This speedup allows for the estimation of increasingly large-scale
connectomes across taxa as these data become available from tracing experiments. The data and code are available online.
\end{abstract}

\textbf{Keywords:} \textit{Matrix equations, computational neuroscience, low-rank approximation, networks}

\section{Introduction}

Neuroscience and machine learning are now enjoying a shared moment of
intense interest and exciting progress.
Many computational neuroscientists find themselves inspired by
unprecedented datasets to develop innovative methods of analysis.
Exciting examples of such next-generation experimental methodology
and datasets are
large-scale recordings and precise manipulations of brain activity,
genetic atlases, and neuronal network tracing efforts.
Thus, techniques which summarize many experiments into an
estimate of the overall brain network are increasingly important.
Many believe that uncovering such network
structures will help us unlock the principles underlying neural computation
and brain disorders~\cite{grillner2016}.
Initial versions of such connectomes \cite{Knox18}
are already being integrated into large-scale modeling projects \cite{reimann2019}.
We present a method which allows us to perform these reconstructions faster,
for larger datasets.

Structural connectivity refers to the synaptic connections
formed between axons (outputs) and dendrites (inputs)
of neurons, which allow them to communicate chemically and electrically.
We represent such networks as a weighted, directed graph
encoded by a nonnegative adjacency matrix $W$.
The network of whole-brain connections or {\it connectome} is
currently studied at a number of scales \cite{sporns2010,kennedy2016}:
Microscopic connectivity catalogues individual neuron connections
but currently is restricted to small volumes due to difficult tracing
of convoluted geometries \cite{kasthuri2015}.
Macroscopic connectivity refers to connections between larger
brain regions and is currently known for a number of model organisms \cite{buckner2019}.
Mesoscopic connectivity
\cite{mitra2014} lies between these two extremes
and captures projection patterns of groups of hundreds to thousands of neurons
among the $10^{6}$--$10^{10}$ neurons in a typical mammalian brain.

Building on previous work \cite{Harris16,Knox18},
we present a scalable method to infer spatially-resolved
mesoscopic connectome from tracing data.
We apply our method to data from the
Allen Mouse Brain Connectivity Atlas \cite{Oh14}
to reconstruct mouse cortical connectivity.
This resource is one of the most comprehensive publicly available
datasets, but similar data are being collected for fly
\cite{jenett2012}, rat \cite{bota2003}, and marmoset \cite{majka2016},
among others.
Our focus is presenting and profiling an improved algorithm for
connectome inference.
By developing scalable methods as in this work,
we hope to enable the reconstruction of high-resolution connectomes
in these diverse organisms.

\subsection{Mathematical formulation of a spatial connectome regression problem}

\begin{figure}[htb!]
  \centering
  \includegraphics[width=\linewidth]{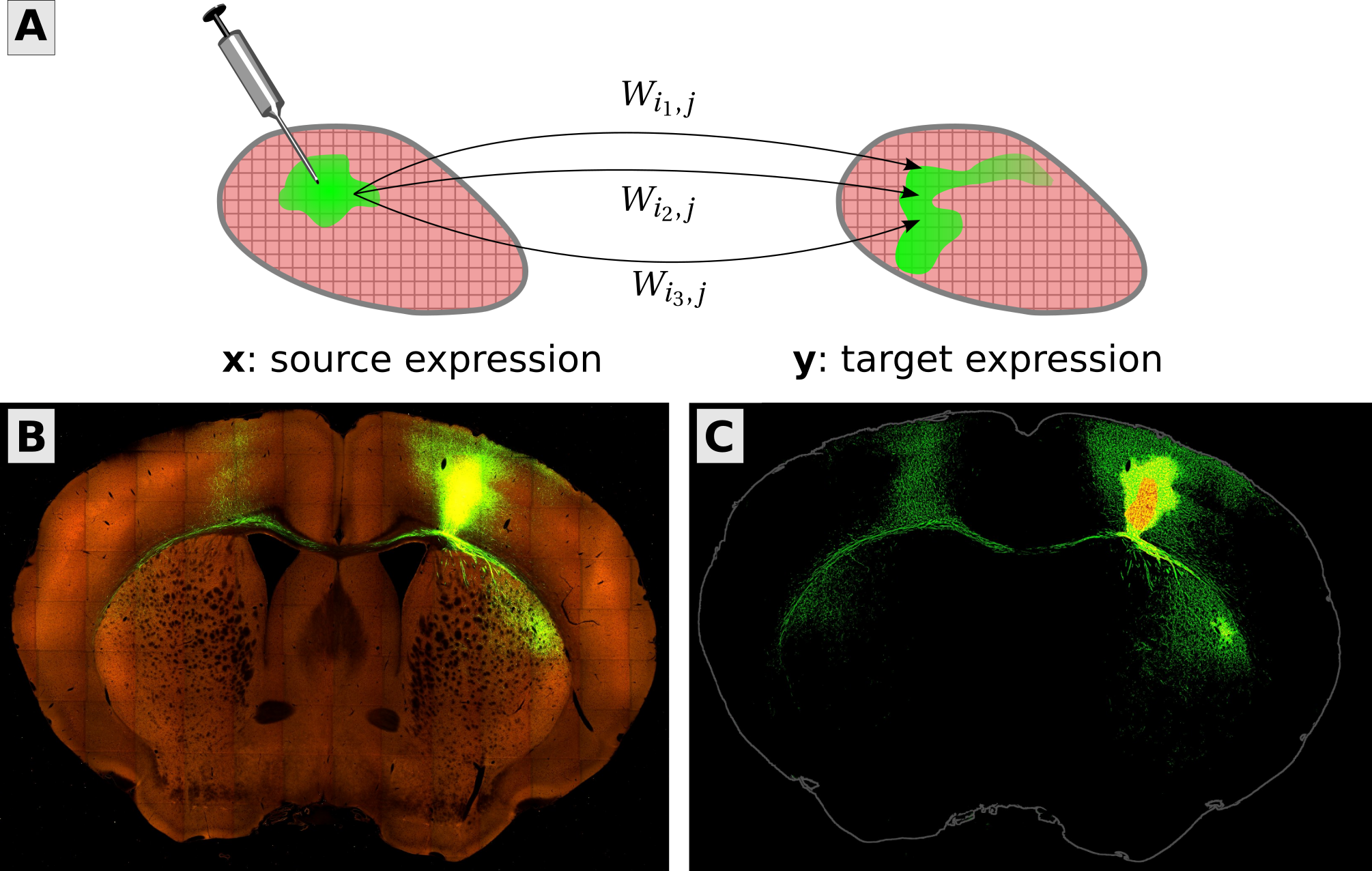}
  \caption{
    In this paper, we present an improved method for
    the mesoscale connectome inference problem.
    (A) The goal is to find a voxel-by-voxel matrix $W$ so that the pattern
    of neural projections $y$ arising from an injection $x$
    is reproduced by matrix-vector multiplication, $y \approx W x$.
    The vectors $x$ and $y$ contain the fraction of fluorescing pixels in each voxel
    from viral tracing experiments.
    {(B)} An example of the data, in this case a coronal slice
    from a tracing experiment delivered to primary motor cortex (MOp).
    Bright green areas are neural cells expressing the green fluorescent protein.
    {(C)} The raw data are preprocessed to separate the injection
    site (red/orange) from its projections (green).
    Fluorescence values in the injection site enter into the source vector $x$,
    whereas fluorescence everywhere else is stored in the entries of
    the target vector $y$.
    The $x$ and $y$ are discretized volume images of the brain
    reshaped into vector form.
    Entry $W_{ij}$ models the expected fluorescence
    at location $i$ from one unit of source fluorescence at location $j$,
    a linear operator mapping from source images to target images.
    Image credit (B and C): Allen Institute for Brain Science.
  }
  \label{fig:data_schematic}
\end{figure}

We focus on the mesoscale because it is naturally captured by
viral tracing experiments
(Figure~\ref{fig:data_schematic}).
In these experiments,
a virus is injected into a specific location in the brain,
where it loads the cells with proteins that can then be imaged,
tracing out the projections of those neurons with
cell bodies located in the injection site.
The source and target signals, within and outside of the injection sites,
are measured as the fraction of fluorescing pixels within cubic voxels.
These form the data matrices
$X\in\R^{\nx\times \ninj}$ and $Y\in\R^{\ny\times \ninj}$,
where parameters $\nx$ and $\ny$ are the number of locations
in the discretized source and target regions of the $d$-D brain,
and $\ninj$ is the number of injections.
In general, $\nx$ and $\ny$ may be unequal, e.g.\ if injections
were only delivered to the right hemisphere of the brain.
Each experiment only traces out the projections from that particular injection site.
By performing many such experiments, with multiple mice,
and varying the injection sites to cover
the brain, one can then ``stitch'' together a mesoscopic connectome
for the average mouse.
We refer the interested reader to \cite{Oh14}
for more details of the experimental procedures.

We present a new low-rank approach to solving
the smoothness-regularized optimization problem
posed by \cite{Harris16}.
Specifically, they considered solving the regularized least squares problem
\begin{equation}
\label{eq:obj_fun_nonneg}
W^* = \arg \min_{W \geq 0} \;
\underbrace{\frac{1}{2} \, \| P_{\Omega} (WX-Y) \|_F^2}_{\text{\small loss}} +
\underbrace{\frac{\lambda}{2} \, \|L_y W + W L_x^\trans\|_F^2}_{\text{\small regularization}},
\end{equation}
where the minimum is taken over nonnegative matrices.
The operator
$P_{\Omega}$
defines an entry-wise product (Hadamard product)
$P_{\Omega} (M) = M \circ \Omega$,
for any matrix $M \in \mathbb{R}^{\ny \times \ninj}$,
and $\Omega$ is a binary matrix, of the same size,
which masks out the injection sites where the entries of $Y$ are
unknown.\footnote{In this paper, we take a different convention for
$\Omega$ (the complement) as in \cite{Harris16}.}
We take the smoothing matrices
$L_y\in\R^{\ny\times \ny}$ and $L_x\in\R^{\nx\times \nx}$
to be discrete Laplace operator,
i.e.\ the graph Laplacians
of the voxel face adjacency graphs
for discretized source and target regions.
We choose a regularization parameter
$\bar\lambda$ and set
$\lambda = \bar\lambda\frac{\ninj}{\nx}$ to avoid dependence on $\nx, \ny$ and $\ninj$,
since the loss term is a sum over $\ny \times \ninj$ entries
and the regularization sums over $\ny \times \nx$ many entries.

We now comment on the typical parameters for problem \eqref{eq:obj_fun_nonneg}.
The mouse brain gridded at 100 \micron resolution
contains approximately
$\nx, \ny
\in \mathcal{O}(10^5)$ voxels in 3-D.
On the other hand, the number of experiments $\ninj$ is less than
$\mathcal{O}(10^3)$.
By projecting the 3-D cortical data into 2-D, as we do in this paper,
we can reduce the size by an order of magnitude to
$\nx, \ny \in\mathcal{O}(10^4)$,
but focusing on the cortex reduces
$\ninj$ to $\mathcal{O}(10^2)$.
Since $\ninj \ll \nx, \ny$,
least squares estimation of $W$
(i.e.\ $\lambda = 0$)
is highly underdetermined
and will remain underdetermined
unless orders of magnitude more tracing
experiments are performed.
Regularization is thus essential for filling the gaps
in injection coverage.
Furthermore, the vast size of the
$\ny \times \nx$ matrix $W$
for whole-brain connectivities has motivated our search for
scalable and fast low-rank methods.

\subsection{Previous methods of mesoscale connectome regression}

Much of the work on mesoscale mouse connectomes
leverages the data and processing pipelines of the
Allen Mouse Brain Connectivity Atlas
available at
\url{http://connectivity.brain-map.org}
\cite{lein2007,Oh14}.
In the early examples of such work,
\cite{Oh14} used viral tracing data to construct regional connectivity
matrices.
Nonnegative matrix regression was used to estimate
the regional connectivity.
First, the injection data were processed into a pair of matrices
$X^\text{Reg}$ and $Y^\text{Reg}$
containing the regionalized injection volumes and projection volumes,
respectively.
The rows of these matrices are the regions and the columns index
injection experiments.
\cite{Oh14} then used nonnegative least squares to fit a
region-by-region matrix $W^\text{Reg}$
such that $Y^\text{Reg} \approx W^\text{Reg} X^\text{Reg}$.
Due to numerical ill-conditioning and a lack of data,
some regions were excluded from the analysis.
Similar techniques have been used to estimate regional connectomes in other animals.
\cite{ypma2016} took a different approach,
using a likelihood-based Markov chain Monte Carlo
method to infer regional connectivity and weight uncertainty from the Allen data.

\cite{Harris16}
made a conceptual and methodological leap when
they presented a method to use such data for spatially-explicit
mesoscopic connectivity.
The Allen Mouse Brain Atlas is essentially a coordinate mapping
which discretizes the average mouse brain into cubic voxels,
where each voxel is assigned to a given region in a hierarchy
of brain regions.
They used an assumption of spatial smoothness to formulate
\eqref{eq:obj_fun_nonneg},
where the specific smoothing term results in a high-dimensional thin plate spline fit
\cite{wahba1990}.
They then solved \eqref{eq:obj_fun_nonneg}
using the generic quasi-Newton algorithm L-BFGS-B \cite{byrd1995}.
This technique was applied to the mouse visual areas
but limited to small datasets since $W$ was dense.
Using a simple low-rank version
based on projected gradient descent,
\cite{Harris16} argued that such a method
could scale to larger brain areas.
However, the initial low-rank implementation
turned out to be too slow to converge for large-scale applications.
Times to convergence were not reported in the original paper,
but the full rank version typically took around a day,
while the low-rank version needed multiple days to reach a
minimum.\footnote{KD Harris, personal communication, 2017.
Note that these times are for the much smaller visual areas dataset.}

\cite{Knox18} simplified the mathematical problem
by assuming that the injections were
delivered to just a single voxel at the injection center of mass.
Using a kernel smoother led to a method which
is explicitly low-rank,
with smoothing performed only in the injection space
(columns of $W$).
This kernel method was applied to the whole mouse brain,
yielding the first estimate of
voxel-voxel whole-brain connectivity for this animal.
However, these assumptions
do not hold in reality:
The injections affect a volume of the brain that encompasses much more
than the center of mass.\footnote{
Wildtype injections can cover 30--500 voxels, approximately 240 on average,
at 100 \micron resolution \cite{Oh14}.}
We also expect that the connectivity is also smooth across projection space
(rows of $W$), since the incoming projections to a voxel are strongly correlated
with those of nearby voxels.
These inaccuracies mean that the kernel method is prone to artifacts,
in particular ones arising from the injection site locations, since there is no ability
for that method to translate the source of projections smoothly away from injection sites.
It is thus imperative to develop an efficient method for the
spline problem that works for large datasets.

\subsection{Continuous formulation motivates the need for sophisticated solvers}

We will now describe, for the first time,
the continuous mathematical properties of this problem,
in order to illuminate why it is challenging to solve.
Equation \eqref{eq:obj_fun_nonneg} can be seen as a discrete
version of an underlying continuous problem (similar to~\cite{rudin1992}, among others)
where we define the cost as
\begin{equation}
\label{eq:costfunctional}
\frac{1}{2}
\sum_{i=1}^\ninj
\int_{T \cap \Omega_i}
\left(
  \int_S \mathcal{W}(x,y) X_i(x) \dee x
  - Y_i(y)
\right)^2 \dee y
+
\frac{\lambda}{2}
\int_T \int_S
\left( \Delta \mathcal{W}(x,y) \right)^2
\dee x \,\dee y .
\end{equation}
The cost is minimized over $\mathcal{W}: T \times S \to \R$,
the continuous connectome,
in an appropriate Sobolev space
(square-integrable derivatives up to fourth order on $T \times S$ is sufficiently regular).
The function $\mathcal{W}$
may be seen as the kernel of an integral operator from $S$ to $T$.
These regions $S$ and $T$ are both compact subsets of $\R^d$
representing source and target regions of the brain.
The mask region $\Omega_i \subset T$ is the subset of the brain
excluding the injection site.
Finally, the discrete Laplacian terms $L$ have been replaced by
the continuous Laplacian operator $\Delta$ on $S \times T$.
The parameter $\lambda$ again sets the level of smoothing.\footnote{
One may consider rescaling $\lambda$ as before,
but subtle differences arise.
In the continuous versus discrete cases
the units of the equations are different,
since the functions $X_i(x)$
and $Y_i(y)$ are now viewed as densities.
Furthermore, there is a mismatch in units between
\eqref{eq:obj_fun_nonneg} and \eqref{eq:costfunctional},
because the graph Laplacian is unitless whereas the Laplace operator is not.
This explains the lack of any dependence on the grid size in the scaling of the discrete problem.
Regardless, choosing the exact scaling to make the continuous and discrete cases match is
not necessary for the more qualitative argument we are making.}

For simplicity, consider
$S = T =$ the whole brain,
$\Omega_i = T$ for all $i = 1, \ldots, \ninj$
and relax the constraint of nonnegativity on $\mathcal{W}$.
Taking the first variational derivative of \eqref{eq:costfunctional}
and setting it to zero yields the Euler-Lagrange equations for
this simplified problem:
\begin{align}
\label{eq:euler-lagrange}
0 &= \lambda
 \Delta^2 \mathcal{W}(x,y)
- \sum_{i=1}^{\ninj} X_i(x) Y_i(y)
+ \int_S \mathcal{W}(x',y)
    \left( \sum_{i=1}^{\ninj} X_i(x') X_i(x) \right) \dee x'
\nonumber
\\
&=
\lambda
 \Delta^2 \mathcal{W}(x,y)
  - g(x,y)
  + \int_S \mathcal{W}(x',y) f(x',x) \dee x' ,
\end{align}
where for convenience we have defined the
data covariance functions
$
  f(x',x) = \sum_{i=1}^{\ninj} X_i(x') X_i(x)
$
and
$
g(x,y) = \sum_{i=1}^{\ninj} X_i(x) Y_i(y),
$
analogous to $X X^T$ and $YX^T$.
The operator $\Delta^2$ is the biharmonic
operator or bi-Laplacian.
Equation \eqref{eq:euler-lagrange} is a fourth-order
partial integro-differential equation in $2d$ dimensions.

Iterative solutions via gradient descent or quasi-Newton methods
to biharmonic and similar equations can be slow to converge
\cite{altas1998}.
It takes many iterations to propagate the highly local action
of the biharmonic differential operator across global spatial scales
due to the small stable step size \cite{rudin1992},
whereas the integral part is inherently nonlocal.
Very slow convergence is what we have found when applying methods
like gradient descent to problem \eqref{eq:obj_fun_nonneg}, also for
low-rank versions.
This includes quasi-Newton methods such as L-BFGS \cite{byrd1995}.
When we attempted to solve the whole-cortex top view and flatmap
problems as in Sections~\ref{sec:test2} and \ref{sec:test3},
the method had not converged (from a naive initialization)
after a week of computation.
These difficulties motivated the development of the method we present here.

\subsection{Outline of the paper}
We present a greedy, low-rank algorithm tailored to the
connectome inference problem.
To leverage powerful linear methods,
we consider solutions to the unconstrained problem
\begin{equation}
\label{eq:obj_fun}
W^* = \arg \min_{W}
\frac{1}{2} \| P_{\Omega} (WX-Y) \|_F^2 +
\frac{\lambda}{2} \|L_y W + W L_x^\trans\|_F^2,
\end{equation}
where all of the matrices and parameters are as in \eqref{eq:obj_fun_nonneg}.
In practice, solutions to the linear problem
\eqref{eq:obj_fun}
are often very close to those of the nonnegative problem
\eqref{eq:obj_fun_nonneg},
since the data matrices $X$ and $Y$ and the ``true'' $W$ are nonnegative.
Setting any negative entries in the computed solution
$W^*$ to zero is adequate,
or can serve as an initial guess
to an iterative solver for the slower nonnegative problem.

Equation \eqref{eq:obj_fun} is another regularized
least-squares problem.
In Section~\ref{sec:gradient_low_rank},
we show that taking the gradient
and setting it equal to zero leads to a linear matrix equation
in the unknown $W$.
This takes the form of a generalized Sylvester equation with
coefficient matrices formed from the data and Laplacian terms.
The data matrices are, in fact, of low-rank since $\ninj \ll \nx, \ny$,
and thus we can expect a low-rank approximation
$W \approx UV^\trans$ to the full solution to perform well (see~\cite{Harris16}, although we do not know how to justify this rigorously).
We provide a brief survey of some low-rank methods for linear
matrix equations in Section~\ref{sec:low_rank_methods}.
We employ a greedy solver that finds rank-one components $u_i v_i^\trans$
one at a time, explained in Section~\ref{sec:greedy}.
After a new component is found, it is orthogonalized and a
Galerkin refinement step is applied.
This leads to Algorithm~\ref{alg:greedylr}, our complete method.

We then test the method on a few connectome fitting problems.
First, in Section~\ref{sec:test1},
we test on a fake ``toy'' connectome, where we know the truth.
This is a test problem
consisting of a 1-D brain with smooth connectivity \cite{Harris16}.
We find that the output of our algorithm converges to the true solution
as the rank increases and as the stopping tolerance decreases.
Next, we present two benchmarks using real viral tracing data
from the isocortices of mice,
provided by the Allen Institute for Brain Science.
In each case, we work with 2-D data in order to limit the problem size
and because the cortex is a relatively flat, 2-D shape.
It has also been argued that such a projection also denoises such data
\cite{vanessen2013,gamanut2018}.
In Section~\ref{sec:test2}, we work with data that are averaged directly
over the superior-inferior axis to obtain a flattened cortex.
We refer to this as the {\it top view} projection.
In contrast, for Section~\ref{sec:test3}, the data are flattened
by averaging along curved streamlines of cortical depth.
We call this the {\it flatmap} projection.

Finally, Section~\ref{sec:discussion}
discusses the limitations of our method and
directions for future research.
Our data and code are described in Section~\ref{sec:data}
and freely available for anyone who would like to reproduce the results.

\section{Greedy low-rank method}

\subsection{Linear matrix equation for the unknown connectivity}
\label{sec:gradient_low_rank}
We now derive the equivalent of the ``normal equations''
for our problem.
Denote the objective function \eqref{eq:obj_fun} as $J(W)$,
with decomposition
\[
J(W) = J_\mathrm{loss}(W) + J_\mathrm{reg}(W) =
\frac{1}{2} \| P_{\Omega} (WX-Y) \|_F^2 +
\frac{\lambda}{2} \|L_y W + W L_x^\trans\|_F^2.
\]
Writing $J_\mathrm{loss}$ indexwise,
we obtain (note that $\Omega \circ \Omega=\Omega$)
\[
J_\mathrm{loss} =
\frac{1}{2}\sum_{i,\alpha=1}^{n, \ninj}
\Omega_{i,\alpha}
\left(\sum_{k=1}^{m} W_{i,k} X_{k,\alpha} - Y_{i,\alpha}\right)^2.
\]
The derivative reads
\begin{align*}
\frac{\partial J_\mathrm{loss}}{\partial W_{\hat \imath,\hat k}}
&=
\sum_{i,\alpha=1}^{n,\ninj}
\Omega_{i,\alpha}
\left(
\sum_{k=1}^{m} W_{i,k} X_{k,\alpha} - Y_{i,\alpha}
\right)
X_{\hat k,\alpha} \delta_{i,\hat \imath} \\
&=
\sum_{\alpha=1}^{\ninj}
\Omega_{\hat \imath,\alpha} X_{\hat k,\alpha}
\sum_{k=1}^{m}
\left(
X_{k,\alpha} W_{\hat \imath,k} - X_{\hat k,\alpha}
\Omega_{\hat \imath,\alpha} Y_{\hat \imath,\alpha}
\right) ,
\end{align*}
or in the vector form
$$
\frac{\partial J_\mathrm{loss}}{\partial \vvec (W)}
=
\sum_{\alpha=1}^{\ninj}
\left[
(X_{\alpha} X_{\alpha}^\trans) \otimes \mathrm{diag}(\Omega_{\alpha})
\right]
\vvec (W)
-
\vvec \left((\Omega \circ Y) X^\trans\right),
$$
where $X_{\alpha}$ is the $\alpha$-th column of $X$
and likewise for $\Omega$.
Setting the derivative equal to zeros leads
to the system of normal equations
\begin{equation}
  \label{eq:linear_system}
  A \, \vvec (W) = \vvec \left((\Omega \circ Y) X^\trans\right) ,
\end{equation}
where $\vvec(W)$ is the vector of all columns of $W$ stacked on top of each other.
This linear system features the following
$(\ny \nx) \times (\ny \nx)$ matrix,
consisting of $\ninj+3$ Kronecker products,
\begin{equation}
  \label{eq:system_matrix}
  A =
  \sum_{\alpha=1}^\ninj
  (X_{\alpha} X_{\alpha}^\trans) \otimes \mathrm{diag}(\Omega_{\alpha})
  +
  \lambda (L_x^2 \otimes I_\ny + 2 L_x \otimes L_y + I_\nx \otimes L_y^2) .
\end{equation}
Note that without the observation mask,
$\Omega$ is a matrix of all ones,
and the first term compresses to $XX^\trans \otimes I_\ny$.

The linear system \eqref{eq:linear_system}
can be recast as the linear matrix equation
\begin{align}
\label{multitermSylv}
\cA(W) = D,
\end{align}
with the operator $ \cA(W):=\lambda\cB(W)+\cC(W) $, where
$$
\cB(W):=WL^2_x+2L_yWL_x+L^2_yW,
\quad
\cC(W):=\sum_{\alpha=1}^\ninj
\mathrm{diag} (\Omega_{\alpha}) WX_{\alpha} X_{\alpha}^\trans,
\quad
\mbox{ and }
D:=(\Omega \circ Y) X^\trans.
$$
The smoothing term $\cB$
can be expressed as a squared standard Sylvester operator
$
\cB(W)=\cL(\cL(W)),
$
where
$\cL(W):=L_y W + W L_x.$
The operator $\cL$ is the graph Laplacian operator
on the discretization of $T \times S$.
Furthermore, the right hand side $D$ is a matrix of rank $\ninj$,
since it is an outer product of two rank $\ninj$ matrices.

\subsection{Numerical low-rank methods for linear matrix equations}
\label{sec:low_rank_methods}

Because of the potentially high dimensions $\nx,\ny$,
directly solving the algebraic matrix equation~\eqref{multitermSylv} is numerically inefficient since the
solution will be a dense $\ny\times \nx$ matrix, making even storing it infeasible.
However,
the rank of the right hand side of \eqref{multitermSylv}
is at most $\ninj \ll \nx,\ny$.
It is often observed and theoretically shown~\cite{Gra04,BenB13,JarMPR17}
that the solutions of large matrix equations with low-rank right hand
sides exhibit rapidly decaying singular values.
Hence, the solution $W$ is expected to have small
numerical rank in the sense that few of its
singular values are larger than machine precision or the experimental noise floor.
Intuitively, since we also seek very smooth solutions,
this also helps control the rank, since high frequency components
tend to be associated with small singular values.
This motivates us to
approximate the solution of \eqref{multitermSylv} by a low-rank approximation
$W\approx UV^\trans$
with
$U\in\R^{\ny \times \rank}$,
$V\in\R^{\nx\times \rank}$
and $\rank \ll \min(\nx,\ny)$.
The low-rank factors are then typically
computed by iterative methods which never form the approximation $UV^\trans$ explicitly.

Several low-rank methods for computing $U,V$ have been proposed, starting from methods for
standard Sylvester equations $AX + XB = D$, e.g.,
\cite{Ben04,BenLT09,BenS13,Sim16}
and more recently
for general linear matrix equations like \eqref{multitermSylv}
\cite{Dam08,BenB13,ShaSS15,JarKMR16,JarMPR17,PowSS17}.
However, these methods are specialized
and require the problem to have particular
structures or properties
(e.g., $\cB,~\cC$ have to form a convergent splitting of $\cA$),
which are not present in the problem at hand.
The main structures present in \eqref{multitermSylv} are positive definiteness
and sparsity of $L_x, L_y$.

An approach that
is applicable to the matrix equation~\eqref{multitermSylv}
is a greedy method
as proposed by \cite{KreS15},
which is based on successive rank-1 approximations of the error.
Because this method is quite general,
we tailored specific components of the algorithm to our problem.
Three main challenges were overcome:
First, we choose a simpler stopping criterion for the ALS routine.
Second, specific solvers were chosen for the
three main sub-problems of the algorithm,
which maximizes its efficiency.
Third, we developed a GPU implementation of the Galerkin refinement,
to make this bottleneck step more efficient.
We advocate for this method in the rest of the paper.
\subsection{Description and application of the greedy low-rank solver}
\label{sec:greedy}
Here we briefly review the algorithm from~\cite{KreS15}
and explain how it is specialized for our particular problem.
Assume there is already an approximate solution
$W_j\approx W^*$
of the linear matrix equation
$\cA(W) = D$,
equation~\eqref{multitermSylv},
with solution $W^*$.
We will improve our solution by an update of rank one:
$W_{j+1}=W_j+u_{j+1}v_{j+1}^\trans$,
where
$u_{j+1}\in\R^{\ny}$ and
$v_{j+1}\in\R^{\nx}$.
The update vectors
$u_{j+1},~v_{j+1}$
are computed by minimizing an error functional that we will soon define.
Since
the operator $\cA$
is positive definite,
it induces the $\cA$-inner product
$\langle X,Y\rangle_{\cA} =
\trace{Y^{\trans}\cA(X)}$
and the $\cA$-norm $\|Y\|_{\cA}:=\sqrt{\langle Y,Y\rangle_{\cA}}$.
So we find $u_{j+1},v_{j+1}$
by minimizing the squared error in the $\cA$-norm:
\begin{align*}
  \left(u_{j+1},v_{j+1}\right) &=
    \arg \min_{u,v} \|W^* - W_j - uv^\trans \|_{\cA}^2
\\
  &= \arg \min_{u,v} \trace{ (W^* - W_j - uv^\trans)^\trans \cA (W^* - W_j - uv^\trans) }
\nonumber
\\
  &= \arg \min_{u,v}
    \trace{ (W^* - W_j - uv^\trans)^\trans (D - \cA (W_j) - \cA(uv^\trans)) }.
\nonumber
\end{align*}
Discarding constant terms,
noting that
$\langle X, Y \rangle_\cA = \langle Y,X \rangle_\cA$,
and setting
$R_j=D-\cA(W_j)$ leads to
\begin{align}
\label{minuv2}
 \left(u_{j+1},v_{j+1}\right) =
    \arg \min_{u,v} \langle uv^\trans,uv^\trans\rangle_{\cA}-2\trace{uv^\trans R_j}.
\end{align}
Notice that the rank-1 decomposition $u v^\trans$ is not unique,
because
we can rescale the factors by any nonzero scalar $c$
such that $(uc)(v/c)^\trans$ represents the same matrix.
This reflects the fact that the
optimization problem~\eqref{minuv2} is not convex.
However, it is convex in each of the factors $u$ and $v$ separately.

We obtain the updates $u_{j+1}$, $v_{j+1}$
via an alternating linear (ALS) scheme \cite{OrtR00}.
Although we only consider low-rank approximations of matrices here,
ALS methods are also used for computing low-rank approximations of
higher order tensors by means of polyadic decompositions, e.g., 
\cite{harshman70,SorVBDeL13}.
First, a fixed $v$ is used in~\eqref{minuv2}
and a minimizing $u$ is computed which is in the next stage kept fixed and
\eqref{minuv2} is solved for a minimizing $v$.
For a fixed vector $v$ with $\|v\|=1$ the minimizing problem is
\begin{align*}
 \hat u &=
    \arg \min_{u} \langle uv^\trans,uv^\trans\rangle_{\cA}-2\trace{uv^\trans R_j}\\
     \begin{split}
          &=\arg \min_{u} \Big(\lambda\left((u^\trans u)v^\trans L_x^2v+2(u^\trans L_y u)(v^\trans L_xv)+u^\trans
L_y^2u\right)\\
&\qquad\qquad+\sum\limits_{\alpha=1}^\ninj(u^\trans\mathrm{diag}(\Omega_{\alpha})u) (v^\trans X_{\alpha} X_{\alpha}^\trans v)\Big)-2u^\trans R_jv
     \end{split}
     \end{align*}
and, hence, $\hat u$ is obtained by solving the linear system of equations
\begin{subequations}
  \begin{align}
  \label{ALS_linsys}
    \hat A \hat u=R_jv,\quad
    \hat A:= \lambda
    \left((v^\trans L_x^2v)I+2L_y (v^\trans L_xv)+L_y^2 \right)
    +
    \sum_{\alpha=1}^\ninj
    \mathrm{diag}(\Omega_{\alpha})
    (v^\trans X_{\alpha} X_{\alpha}^\trans v).
  \end{align}
  The second half iteration starts from the fixed
  $u=\hat u/\|\hat u\|$
  and tries to find a minimizing $\hat v$ by solving
  \begin{align}
    \label{ALS_linsys_v}
    \hat B \hat v=R^\trans_ju,\quad
    \hat B:=\lambda
    \left( L_x^2+2L_x (u^\trans L_yu)+(u^\trans L_y^2u)I \right)
    +
    \sum_{\alpha=1}^\ninj
    (u^\trans\mathrm{diag}(\Omega_{\alpha})u)(X_{\alpha} X_{\alpha}^\trans)
  \end{align}
\end{subequations}
which can be derived by similar steps.
The linear systems~\eqref{ALS_linsys} and~\eqref{ALS_linsys_v}
inherit the sparsity from $L_x,L_y$ and $\Omega$.
Therefore they can be solved by sparse direct
or iterative methods.
We use a sparse direct solver for~\eqref{ALS_linsys},
as this was faster than the alternatives.
The coefficient matrix $\hat B$ in~\eqref{ALS_linsys_v}
is the sum of a sparse (Laplacian terms)
and a low-rank (rank $\ninj$ data terms) matrices.
In this case, we solve~\eqref{ALS_linsys_v} using the
Sherman-Morrison-Woodbury formula~\cite{GolV13} and a direct solver for the sparse inversion.

Both half steps form the ALS iteration which should be stopped
when the iterates are close enough to a critical point,
which might be difficult to check.
Here we propose a simpler approach compared to the one in~\cite{KreS15}.
Since we rescale $u$ and $v$ such that $\|u\|_2=\|v\|_2=1$,
the norm of the other factor is equal to the norm of the full matrix.
In other words, $\|\hat u\|_2 = \|\hat uv^\trans\|_2$ after solving for $\hat u$,
and hence $\| \hat u \|_2$ should converge to the norm of the exact solution.
This motivates a simple criterion: we stop the ALS when
$
(1-\delta) \|\hat u\|_2 \le \|\hat v\|_2 \le (1+\delta) \|\hat u\|_2,
$
where $\hat u$ and $\hat v$ are taken from two consecutive ALS steps,
and $\delta<1$ is a small threshold.
It turns out that a relatively crude tolerance of $\delta=0.1$,
corresponding to 2--4 ALS iterations,
is sufficient in practice for the overall convergence of the algorithm.

The second stage of the method is a non-greedy
Galerkin refinement of the low-rank factors.
Suppose a rank $j$ approximation
$W_j=\sum_{i=1}^j u_i v_i^\trans$ of $W$ has been already computed.
Let $U \in \R^{\ny\times j}$ and $V\in\R^{\nx\times j}$
have orthonormal columns,
spanning the spaces
$\mathrm{span}\lbrace u_1,\ldots,u_j\rbrace$
and
$\mathrm{span}\lbrace v_1,\ldots,v_j\rbrace$,
respectively.
We compute a refined approximation $UZV^\trans$ for $Z\in\R^{j\times j}$
by imposing the following condition onto the residual:
\begin{align*}
  \mbox{(Galerkin condition)\qquad Find $Z$ so that} \quad
  \cA(UZV^\trans) - D
  \quad
  \perp
  \quad
  \lbrace UZV^\trans\in\R^{\ny\times \nx}, \; Z\in\R^{j\times j}\rbrace.
\end{align*}
This leads to the dense, square matrix equation in $Z$ of dimension $j \leq \rank \ll \nx,\ny$:
\begin{align}\label{Proj_Sylv}
  \lambda \left(Z(V^\trans L^2_xV) + 2(U^\trans L_yU)Z(V^\trans L_x V)+(U^\trans L^2_yU)Z\right) 
  +\sum_{\alpha=1}^\ninj (U^\trans\mathrm{diag}(\Omega_{\alpha}) U)Z(V^\trans X_{\alpha} X_{\alpha}^\trans V)
  = U^\trans D V .
\end{align}

Equation~\eqref{Proj_Sylv}
is a projected version of~\eqref{multitermSylv} and inherits its structure
including the positive definiteness of the operator which acts on $Z$.
Instead of using a direct method to solve~\eqref{Proj_Sylv} as in~\cite{KreS15},
we employ an iterative method similar to~\cite{PowSS17}.
Due to the positive definiteness, the obvious method of choice is a dense,
matrix-valued conjugate gradient method (CG).
Moreover, we reduce the number of iterations significantly
by taking the solution $Z$ from the previous greedy step as an initial guess.
The improved solution $W_{j+1}=UZV^\trans$
yields a new residual $R_{j+1}=D-\cA(W_{j+1})$
onto which the ALS scheme is applied to obtain the next rank-1
updates.
The complete procedure is illustrated in Algorithm~\ref{alg:greedylr}.

\begin{algorithm2e}[t]
  \SetEndCharOfAlgoLine{}
  \SetKwInOut{Input}{Input}
  \SetKwInOut{Output}{Output}
  \caption[Greedy low-rank method]{
    Greedy rank-1 method with Galerkin projection for~\eqref{multitermSylv}}
  \label{alg:greedylr}
  \Input{Matrices $L_x,L_y,X,\Omega,Y$,
    maximal rank $\rank \leq \min(\nx, \ny)$,
    tolerance $0 < \tau \ll 1$}
  \Output{Low-rank approximation of $W = U Z V^\trans$ in factored form}
  Initialize
  $W_0=0$, $R_0=D$, $U_0=V_0=[~]$, $j=0$\;
  \Repeat{$j = \rank$ or $\delta W \le \tau$}{%
    Pick initial vector $v$ for ALS with $\| v \| = 1$,
    then get rank-1 update:\;
    \While{$\delta>0.1$}
    {
      Solve $\hat A\hat u=R_jv$ for $\hat u$
  (sparse direct solver)
      and set $u=\hat u/\|\hat u\|$
      \Comment*[r]{Eqn.~\eqref{ALS_linsys}}
      Solve $\hat B\hat v=R^\trans_j u$ for $\hat v$
  (sparse direct + low-rank update)
      and set $v=\hat v/\|\hat v\|$
      \Comment*[r]{Eqn.~\eqref{ALS_linsys_v}}
      $\delta= \left\vert \frac{\|\hat u\|}{\|\hat v\|} -1 \right\vert$ \;
    }
    $U_{j+1}=\mathrm{orth}([U_{j},u])$,
    $V_{j+1}=\mathrm{orth}([V_{j},v])$
    \Comment*[r]{Orthogonalize new factors}
    Increment rank $j \leftarrow j+1$ \;
    Solve Eqn.~\eqref{Proj_Sylv} for $Z_j$
(CG to tolerance $\tau/2$)
    \Comment*[r]{Galerkin update}
    $R_j = D - \cA( U_j Z_j V_j^\trans )$
    \Comment*[r]{Update residual}
    $\delta W = \| U_j Z_j V_j^\trans  - U_{j-1} Z_{j-1} V_{j-1}^\trans \|_F
    / \| U_j Z_j V_j^\trans \|_F$ \;
  }
\end{algorithm2e}

This Galerkin refinement substantially
improves the greedy approximation, leading to a faster convergence rate \cite{KreS15}.
The ALS stage is primarily used
to sketch the projection bases for the Galerkin solution,
which justifies the limited number of ALS steps.
Use of the Galerkin refinement in low-rank decomposition literature
can be traced back to the greedy approximation in the CP tensor format
\cite{nouy-greedy-2010}, as well as orthogonal matching pursuit
approaches in sparse recovery and compressed sensing \cite{pati-omp-1993}
and deflation strategies in low-rank matrix completion~\cite{HardtW14}.

\section{Performance of the greedy low-rank solver on three problems}

\begin{figure}[t!]
 \centering
 \captionof{table}{
   Computing times and errors for the toy brain test problem.
   The output $W$ is compared to truth and a rank 140 reference solution.
 }
 \label{tab:test}
\begin{tabular}{c|ccccc}
 $\mathrm{rank}(W)$          & 10        & 20        & 40        & 60          & 80       \\\hline  
 CPU time (sec.)             & 0.0396    & 0.1554    & 0.9653    & 2.6398      & 3.1108   \\\hline
 $\mathcal{E}(W,W_{140})$    & 3.2324e-01    & 5.5407e-02    & 1.4162e-02    & 1.2125e-03      & 3.1549e-04  \\
 $\mathcal{E}(W,W_{\rm true})$   & 2.9418e-01    & 7.9921e-02    & 7.1537e-02    & 6.9777e-02      & 6.9821e-02  \\\hline
 $\mathcal{E}_{\rm rel}(W, W_{140})$    & 4.3320e-01 &  8.9700e-02 &  2.4900e-02  & 2.5000e-03 &  5.1300e-04 \\
 $\mathcal{E}_{\rm rel}(W,W_{\rm true})$   & 4.0130e-01 &  1.1410e-01 &  1.0350e-01  & 1.0040e-01 &  1.0040e-01 \\
\end{tabular}

\vspace{1em}
 \captionof{figure}{
   Toy brain test problem.
   Top: true connectivity map $W_{\rm true}$ (left)
   and the low-rank solution with $\mathrm{rank}=40$ and $\lambda=100$ (right).
   Bottom: solutions with $\Omega=1$ (left) and $\lambda=0$ (right).
   The locations of simulated injections are shown by the grey bars.
   This shows that both the mask ($\Omega$) and smoothing ($\lambda > 0$)
   are necessary for good recovery.
 }
 \label{fig:test}
 \includegraphics[width=0.45\linewidth, trim=14 3 14 2, clip]
 {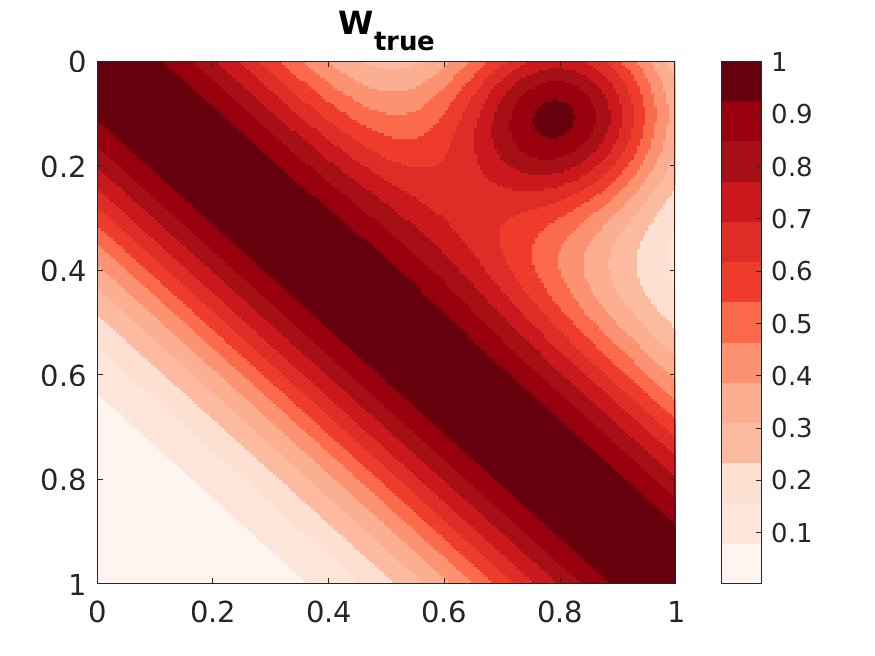} \hskip 3mm
 \includegraphics[width=0.45\linewidth, trim=14 3 14 2 , clip]
 {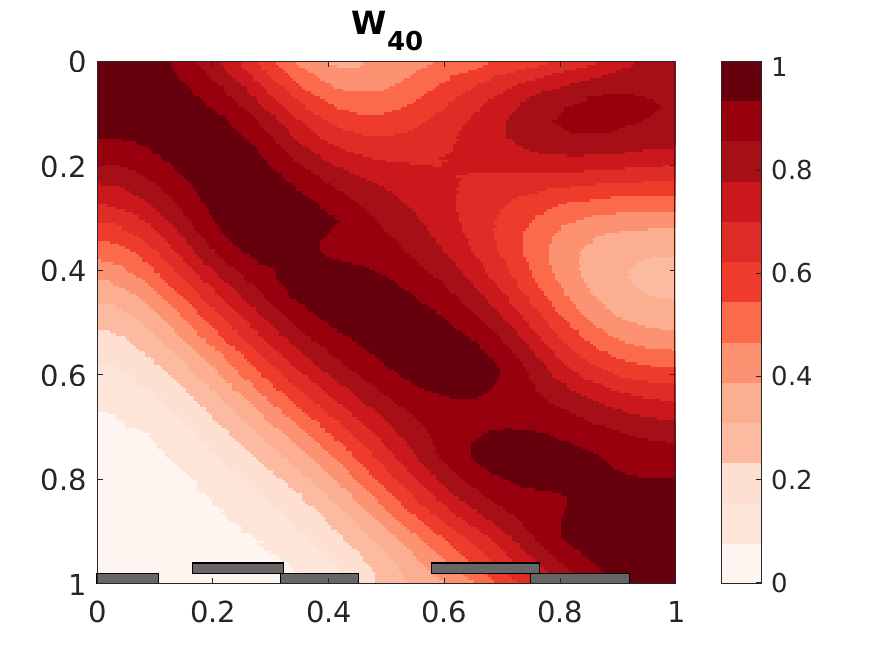} \\
 \vspace{2mm}
 \includegraphics[width=0.45\linewidth, trim=14 3 14 2 , clip]
 {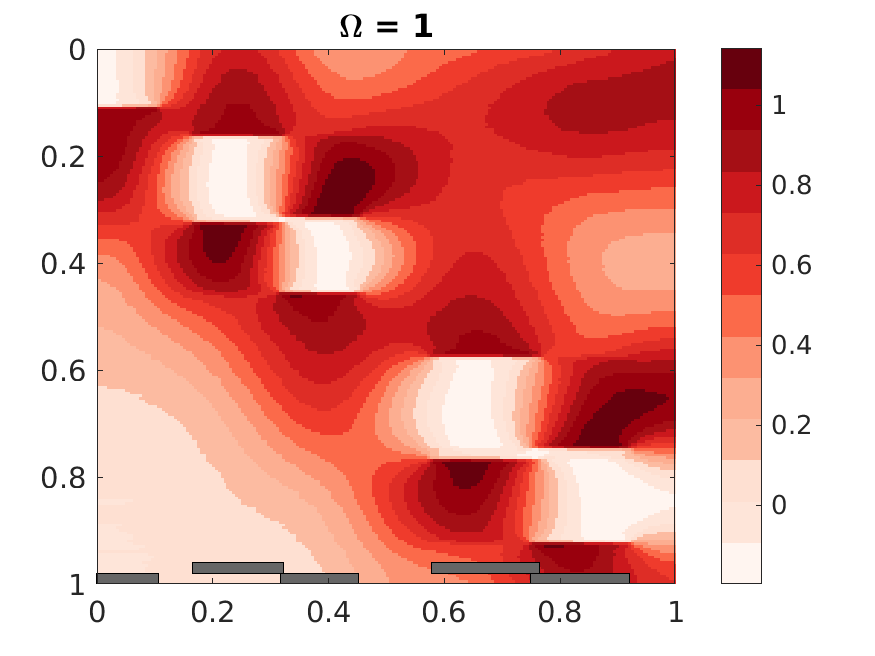} \hskip 3mm
 \includegraphics[width=0.45\linewidth, trim=14 3 14 2 , clip]
 {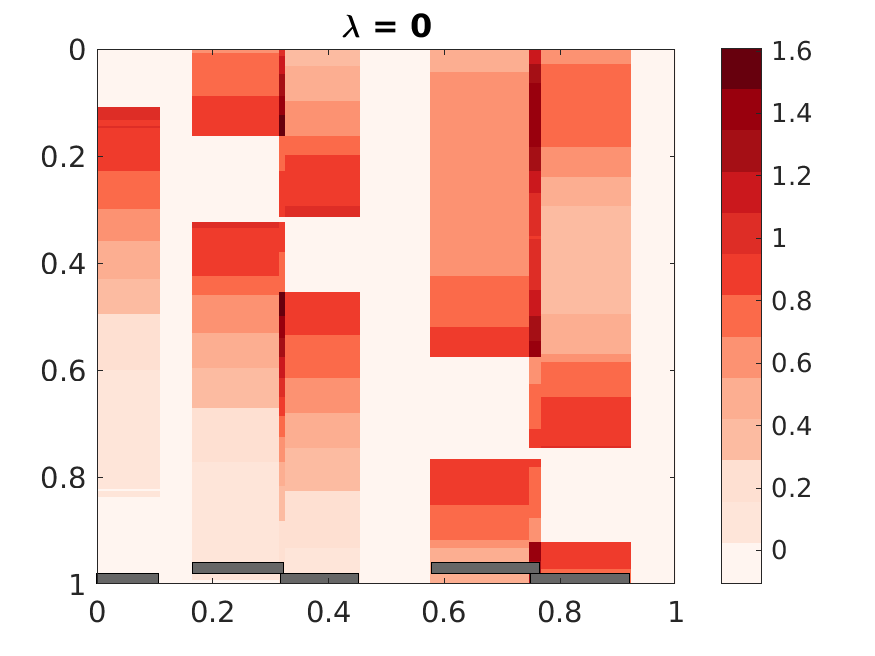}
\end{figure}

\label{sec:tests}

There are three test problems to which we apply Algorithm~\ref{alg:greedylr}:
a toy problem with synthetic data (Section~\ref{sec:test1}),
the top view projected mouse connectivity data (Section~\ref{sec:test2}),
and the flatmap projected data (Section~\ref{sec:test3}).
These tests show that the method easily scales to whole-brain connectome
reconstruction.

We investigate the computational complexity
and convergence of the greedy algorithm.
Since the matrices in~\eqref{ALS_linsys} are sparse, the ALS steps need
$\mathcal{O}(n\rank^2 \ninj)$
operations in total for the final solution rank $\rank$,
where $n = \max(\nx,\ny)$.
In turn, if the solution of~\eqref{Proj_Sylv} takes $\gamma$
CG iterations, this step will have a cost of
$\mathcal{O}(\gamma \rank^3 \ninj)$.
Although $\gamma$ can be kept at the same level for all $j$,
it depends on the stopping tolerance $\tau$, as does the rank $\rank$.
We will therefore investigate the cost
in terms of the total computation time and the
corresponding solution accuracy for a range of solution rank values.

The numerical experiments were performed on an \intel~E5-2650~v2~CPU with
$8$ threads and 64Gb~RAM.
We employ an \nvidia~P100~GPU card for some subtasks:
The Galerkin update relies on
dense linear algebra to solve~\eqref{Proj_Sylv}
by the CG method,
so this stage admits an efficient GPU implementation.
Algorithm~\ref{alg:greedylr} is implemented in
\matlab~R2017b, and was
run on the Balena High Performance Computing Service
at the University of Bath.
See Section~\ref{sec:data} for additional data and code resources.

We measure errors in the solution
using the root mean squared error.
Given any reference solution
$W_\star$ of size $\ny \times \nx$,
e.g.\ the truth or a large-rank solution when the truth is unknown,
and a low-rank solution $W_\rank$,
the RMS error is computed as
$
\mathcal{E}(W_r,W_\star) =
\frac{\left\|W_\rank - W_\star\right\|_F }{\sqrt{\ny \nx}}.
$
We also report the relative error in the Frobenius norm
$
\mathcal{E}_{\rm rel}(W_\rank,W_\star) =
\frac{\left\|W_\rank - W_\star\right\|_F }{\left\|W_\star\right\|_F} .
$

\subsection{Test problem: a toy brain}

\label{sec:test1}

We use the same test problem as in \cite{Harris16},
a one-dimensional ``toy brain.''
The source and target space are $S = T = [0,1]$.
The true connectivity kernel corresponds to a Gaussian
profile about the diagonal
plus an off-diagonal bump:
\begin{equation}
\label{eq:W_true_toy}
W_{\rm true} (x,y)
=
e^{-\left(\frac{x-y}{0.4}\right)^2}
+
0.9\,
e^{ -\frac{\left(x-0.8\right)^2 + (y-0.1)^2}{(0.2)^2} }
 \;.
\end{equation}
The input and output spaces
were discretized using
$\nx = \ny = 200$ uniformly lattice points.
Injections are delivered at $\ninj = 5$
locations in $S$, with a width of
$0.12 + 0.1 \epsilon$,
where
$\epsilon \sim \mathrm{Uniform}(0,1)$.
The values of $X$ are set to 1 within the injection region and
0 elsewhere, $\Omega_{ij} = 1 - X_{ij}$,
$Y$ is set to 0 within the injection region,
and we add Gaussian noise with standard deviation $\sigma = 0.1$.
The matrices $L_x = L_y$ are the 3-point graph Laplacians for the 1-d chain.

We depict the true toy connectivity
$W_\text{true}$ as well as a number of low-rank solutions
output by our method in Figure~\ref{fig:test}.
Both the mask and regularization are required for good performance:
If we remove the mask, setting $\Omega$ equal to the matrix of all ones,
then there are holes in the data at the location of the injections.
If we try fitting with $\lambda = 0$, i.e.\ no smoothing,
then the method cannot fill in holes or extrapolate outside the injection sites.
It is only with the combination of all ingredients
that we recover the true connectivity.

In Table~\ref{tab:test} we show the performance of the algorithm for
ranks $r = 10$, 20, 40, 60, and 80.
The output $W$ is compared to $W_\text{true}$ as well as
the rank 140 output of the algorithm.
The stopping tolerance was $\tau = 10^{-7}$ to ensure that the algorithm has reached this maximal rank.
We see that the RMS distance to the reference solution $W_{140}$
decreases as we increase the rank,
as does the relative distance.
However, the RMS and relative distances from $W_\text{true}$
asymptote to roughly $0.07$ and 10\%, respectively, by rank 40.
This shows that rank 40 is a suitable maximum rank for this problem given the
input data and noise.

The computing time of the greedy method (in this example we use the CPU only version)
remains in the order of seconds even
for the largest considered ranks.
In contrast, the unpreconditioned CG method needs thousands of iterations
(and hundreds of second of time) to compute a solution within the same order of accuracy.
Since it is unclear how to develop a preconditioner for Eq. \eqref{eq:linear_system},
especially for a non-trivial $\Omega$,
in the next sections we focus only on the greedy algorithm.

\subsection{Mouse cortex: top view connectivity}

\label{sec:test2}

\begin{figure}[tp!]
 \centering
  \captionof{figure}{
    Computing times and errors for the top view data.
    The errors are computed with reference to the full-rank
    solution $W_\star = W_{\rm full}$.
    Full rank time: $\gg 6 \times 10^5$ s (see text).
  }
 \label{fig:top_view_times}
 \resizebox{0.49\linewidth}{!}{%
 \begin{tikzpicture}
  \begin{axis}[%
  xmode=normal,
  ymode=log,
  xmin=100,xmax=1000,
  xlabel=$r$,
  ylabel={CPU time (seconds) of low-rank method},
  legend style={at={(0.99,0.01)},anchor=south east},
  ]

  \addplot+[] coordinates{
   (125 , 5.9397e+02)
   (200 , 1.3559e+03)
   (250 , 2.1921e+03)
   (400 , 5.5883e+03)
   (500 , 1.0658e+04)
   (750 , 5.8047e+04)
   (1000, 2.2193e+05)
  }; \addlegendentry{CPU};

  \addplot+[] coordinates{
   (125 , 4.0837e+02)
   (200 , 8.2131e+02)
   (250 , 9.2984e+02)
   (400 , 1.8375e+03)
   (500 , 2.3226e+03)
   (750 , 5.8386e+03)
   (1000, 1.3494e+04)
  }; \addlegendentry{GPU};
  \end{axis}
 \end{tikzpicture}}~
 \resizebox{0.49\linewidth}{!}{%
 \begin{tikzpicture}
  \begin{axis}[%
  xmode=normal,
  ymode=log,
  xmin=100,xmax=1000,
  xlabel=$r$,
  ylabel={Error of low-rank method},
  legend style={at={(0.99,0.99)},anchor=north east},
  ]

  \addplot+[] coordinates{
   (125 , 7.9451e-05)
   (200 , 2.3655e-05)
   (250 , 1.1267e-05)
   (400 , 1.8239e-06)
   (500 , 7.2494e-07)
   (750 , 8.1841e-08)
   (1000, 4.4854e-09)
  }; \addlegendentry{$\mathcal{E}(W_r, W_{\star})$};

  \addplot+[] coordinates{
   (125 , 6.9712e-02)
   (200 , 2.0756e-02)
   (250 , 9.8863e-03)
   (400 , 1.6004e-03)
   (500 , 6.3608e-04)
   (750 , 7.1810e-05)
   (1000, 3.9356e-06)
  }; \addlegendentry{$\mathcal{E}_{\rm rel}(W_r, W_{\star})$};
  \end{axis}
 \end{tikzpicture}}
\end{figure}
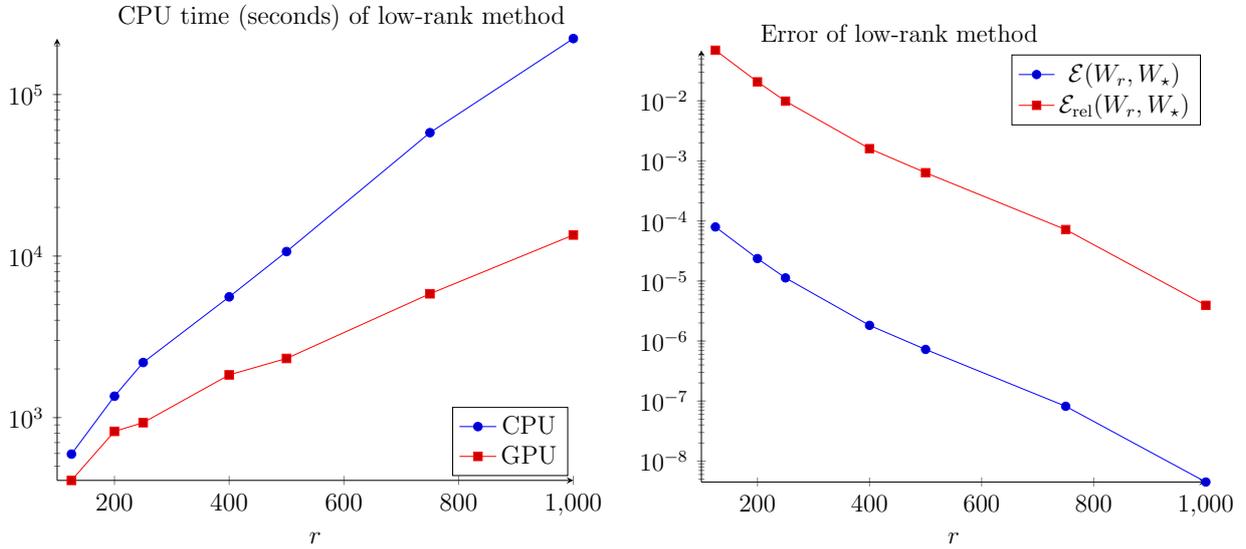


\begin{figure}[tp!]
 \centering
  \captionof{figure}{
    Computing times and errors for the flatmap data.
    The errors are computed with reference to the
    rank-1000 solution $W_\star = W_{1000}$.
  }
 \label{fig:flatmap_times}
 \resizebox{0.49\linewidth}{!}{%
 \begin{tikzpicture}
  \begin{axis}[%
  xmode=normal,
  ymode=log,
  xmin=100,xmax=1000,
  xlabel=$r$,
  ylabel={CPU time (seconds) of low-rank method},
  legend style={at={(0.99,0.01)},anchor=south east},
  ]

  \addplot+[] coordinates{
   (125 , 1.5330e+03)
   (200 , 3.5203e+03)
   (250 , 5.2635e+03)
   (400 , 1.3293e+04)
   (500 , 2.2939e+04)
   (750 , 1.0025e+05)
   (1000, 4.0232e+05)
  }; \addlegendentry{CPU};

  \addplot+[] coordinates{
   (125 , 9.2371e+02)
   (200 , 1.4444e+03)
   (250 , 2.1545e+03)
   (400 , 4.0042e+03)
   (500 , 5.4159e+03)
   (750 , 1.1593e+04)
   (1000, 2.6473e+04)
  }; \addlegendentry{GPU};
  \end{axis}
 \end{tikzpicture}}~
 \resizebox{0.49\linewidth}{!}{%
 \begin{tikzpicture}
  \begin{axis}[%
  xmode=normal,
  ymode=log,
  xmin=100,xmax=1000,
  xlabel=$r$,
  ylabel={Error of low-rank method},
  legend style={at={(0.99,0.99)},anchor=north east},
  ]

  \addplot+[] coordinates{
   (125 , 3.8057e-05)
   (200 , 2.9520e-05)
   (250 , 7.0891e-06)
   (400 , 1.8524e-06)
   (500 , 8.6443e-07)
   (750 , 1.1149e-07)
  }; \addlegendentry{$\mathcal{E}(W_r, W_{\star})$};

  \addplot+[] coordinates{
   (125 , 8.8936e-02)
   (200 , 6.8986e-02)
   (250 , 1.6567e-02)
   (400 , 4.3289e-03)
   (500 , 2.0201e-03)
   (750 , 2.6053e-04)
  }; \addlegendentry{$\mathcal{E}_{\rm rel}(W_r, W_{\star})$};
  \end{axis}
 \end{tikzpicture}}
\end{figure}
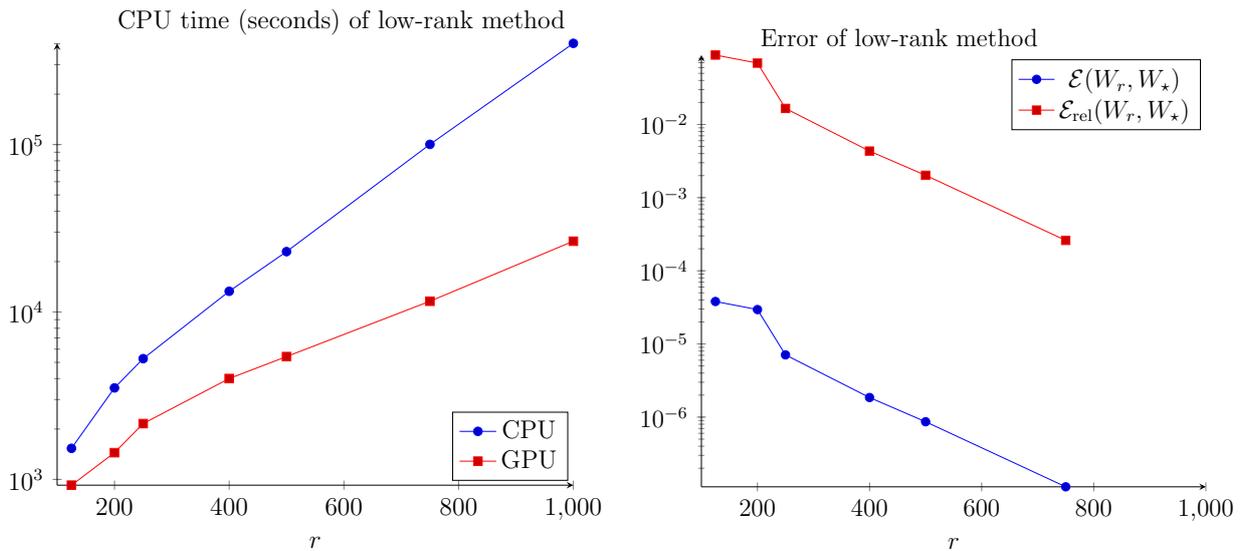


\begin{figure}[tp!]
 \centering
  \captionof{figure}{
    Value of cost function $J(W_r)$ versus the rank $r$ of the low-rank approximation $W_r$.
  }
 \label{fig:J}
 \resizebox{0.49\linewidth}{!}{%
 \begin{tikzpicture}
  \begin{axis}[%
  xmode=normal,
  ymode=normal,
  ymin=4600,ymax=4700,
  xmin=100,xmax=1000,
  xlabel=$r$,
  ylabel={$J(W_r)$},
  legend style={at={(0.99,0.99)},anchor=north east},
  ]

  \addplot+[] coordinates{
   (125 , 4.6740e+03)
   (200 , 4.6410e+03)
   (250 , 4.6376e+03)
   (400 , 4.6361e+03)
   (500 , 4.6360e+03)
   (750 , 4.6360e+03)
   (1000, 4.6360e+03)
  }; \addlegendentry{top view};
  \end{axis}
 \end{tikzpicture}}~
 \resizebox{0.49\linewidth}{!}{%
 \begin{tikzpicture}
  \begin{axis}[%
  xmode=normal,
  ymode=normal,
  ymin=10000,ymax=10400,
  xmin=100,xmax=1000,
  xlabel=$r$,
  ylabel={$J(W_r)$},
  legend style={at={(0.99,0.99)},anchor=north east},
    y tick label style={
        /pgf/number format/.cd,
            fixed,
            scaled y ticks = false,
            precision=4,
        /tikz/.cd
    },
  ]

  \addplot+[] coordinates{
   (125 , 1.0273e+04)
   (200 , 1.0112e+04)
   (250 , 1.0088e+04)
   (400 , 1.0077e+04)
   (500 , 1.0076e+04)
   (750 , 1.0076e+04)
   (1000, 1.0076e+04)
  }; \addlegendentry{flatmap};
  \end{axis}
 \end{tikzpicture}}
\end{figure}
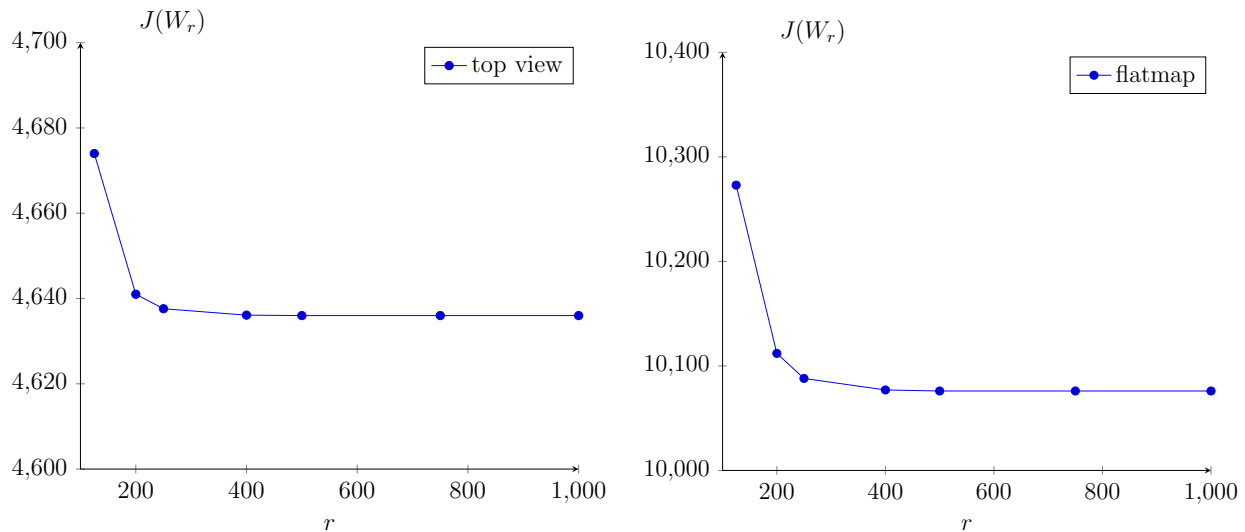

\begin{figure}[t!]
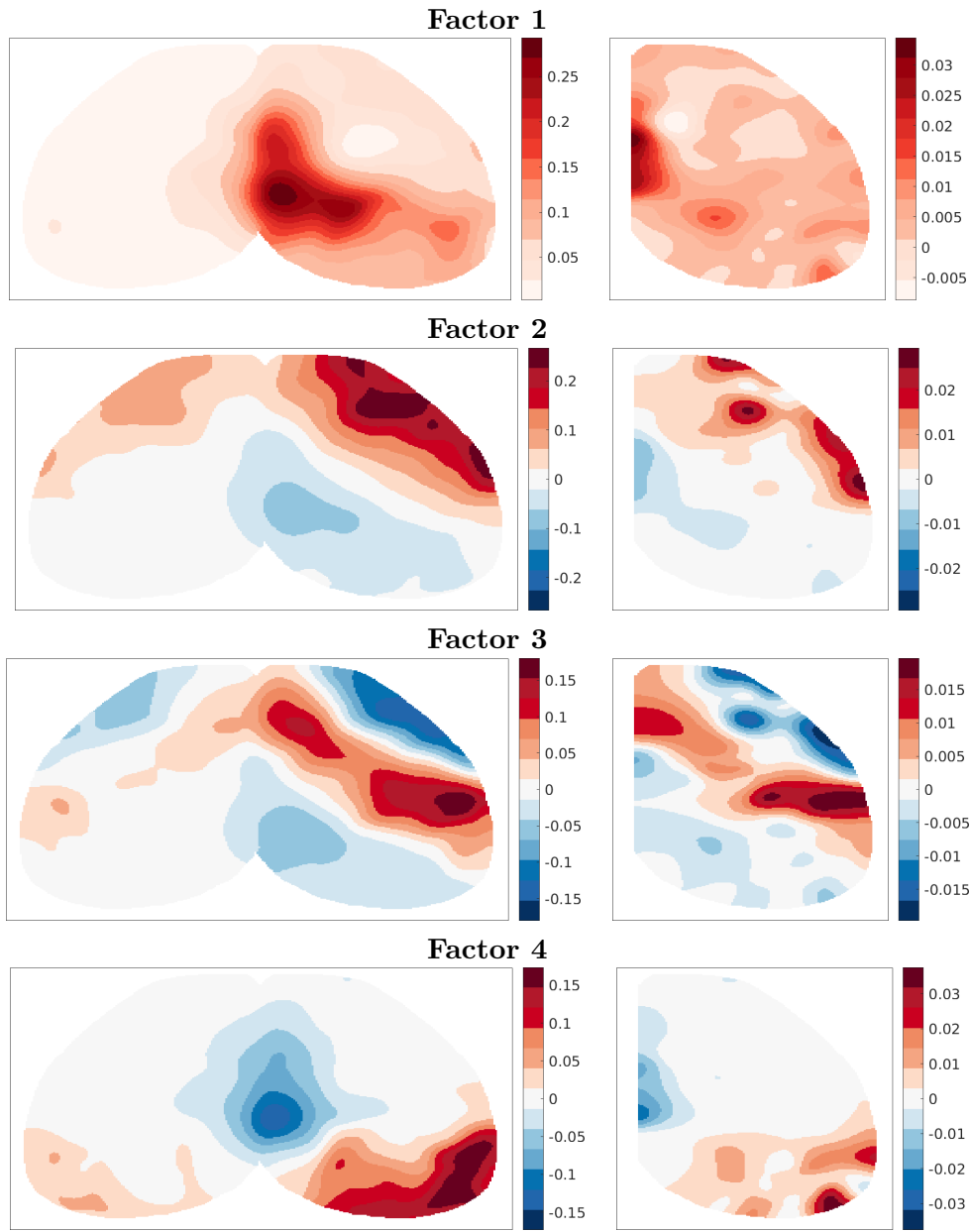

  \centering

  \captionof{figure}{
    Top four singular vectors of the
    top view connectivity with $\rank = 500$.
    Left: scaled target vectors  $\hat U \Sigma$.
    Right: source vectors $\hat V$.}
  \label{fig:top_view}
  \foreach \i in {1,...,4}{
    {\small \bf Factor \i} \\
    \vskip 1mm
    \includegraphics[height=3.5cm]{top_view_US_\i.png}
    \hskip 3mm
    \includegraphics[height=3.5cm]{top_view_V_\i.png}
    \\
  }
\end{figure}

We next benchmark Algorithm~\ref{alg:greedylr} on
mouse cortical data projected into a top-down view.
See Section~\ref{sec:data} for details about how we obtained these data.
Here, the problem sizes are
$\ny = 44\,478$ and $\nx = 22\,377$
and the number of injections $\ninj=126$.
We use the smoothing parameter
$\bar\lambda=10^6$.

We run the low-rank solver with the target solution rank varying from $\rank = 125$ to $1000$.
The stopping tolerances $\tau$
were decreased geometrically from $10^{-3}$ for $r=125$ to $10^{-6}$ for $r=1000$.
This delivers accurate but cheap solution to the Galerkin system \eqref{Proj_Sylv}
while ensuring that the algorithm reached the target rank.

These low-rank solutions are compared to the full-rank solution
$W_\text{full}$
with $\rank = \nx = 22\,377$
found by L-BFGS
\cite{byrd1995},
similar to
\cite{Harris16},
which used L-BFGS-B to deal with the nonnegativity constraint.
Note that this full rank algorithm
was initialized from the output of the low-rank algorithm.
This led to a significant speedup:
The full rank method, initialized naively, had not reached a similar value
of the cost function~\eqref{eq:obj_fun} after a {\it week} of computation,
but this ``warm start'' allowed it to finish within hours.

The computing times and errors are presented in Figure~\ref{fig:top_view_times}.
We see that the RMS errors are relatively small for ranks above 500,
below $10^{-6}$.
Neither the RMS or relative
error seem to have plateaued at rank 1000, but they are small.
At rank 1000, the vector $\ell_\infty$ error
(maximum absolute deviation of the matrices as vectors,
not the matrix $\infty$-norm)
$\| W_{1000} - W_\text{full} \|_\infty$ is less than $10^{-6}$,
which is certainly within experimental uncertainty.
In Figure~\ref{fig:J}, the value of the cost function $J(W_r)$ is plotted against the rank $r$ of the approximation $W_r$ for the top view and flatmap data. Apparently, around $r=500$ the cost function begins to stagnate indicating that the approximation quality does not significantly improve any more. Hence, we continue the investigation with the numerical rank set to $r=500$.

We analyze the leading singular vectors of the solution.
The output of the algorithm is
$W_r = U Z V^\trans$, which is {\em not} the SVD of $W_r$ because $Z$ is not diagonal.
We perform a final SVD of the Galerkin matrix,
$Z = \tilde U \Sigma \tilde V^\trans$ and set
$\hat U = U \tilde U$
and
$\hat V = V \tilde V$,
so that
$W_r = \hat U \Sigma \hat V^\trans$
is the SVD of the solution.

The first four of these singular vectors are shown in
Figure~\ref{fig:top_view}.
The brain is oriented with the medial-lateral axis aligned
left-right and anterior-posterior axis aligned top-bottom,
as in a transverse slice.
The midline of the cortex is in the center of the target plots,
whereas it is on the left edge of the source plots.
We observe that the leading component is a strong projection
from medial areas of the cortex near the midline to nearby locations.
The second component provides a correction which adds
local connectivity among posterior areas and anterior areas.
Note that increased anterior connectivity arises from negative entries
in both source and target vectors.
The sign change along the roughly anterior-posterior axis
manifests as a reduction in connectivity from anterior to posterior regions
as well as from posterior to anterior regions.
The third component is a strong local connectivity among
somatomotor areas located medially along the anterior-posterior axis
and stronger on the lateral side
where the barrel fields, important sensory areas for whisking, are located.
Finally, the fourth component is concentrated in posterior locations,
mostly corresponding to the visual areas,
as well as more anterior and medial locations
in the retrosplenial cortex (thought to be a memory and association area).
The visual and retrosplenial parts of the component show opposite signs,
reflecting stronger local connectivity within these regions than distal
connectivity between them.

These patterns in Figure~\ref{fig:top_view} are reasonable,
since connectivity in the brain is
dominantly local with some specific long-range projections.
We also observe that the projection patterns
(left components $\hat U \Sigma$) are fairly symmetric across the midline.
This is also expected due to the mirroring of major brain areas
in both hemispheres,
despite the evidence for some lateralization, especially in humans.
The more specific projections between brain regions will show up in
later, higher frequency components.
However, it becomes increasingly difficult to
interpret lower energy components as specific pathways,
since these combine in complicated ways.

\subsection{Mouse cortex: flatmap connectivity}

\label{sec:test3}

\begin{figure}[tp!]
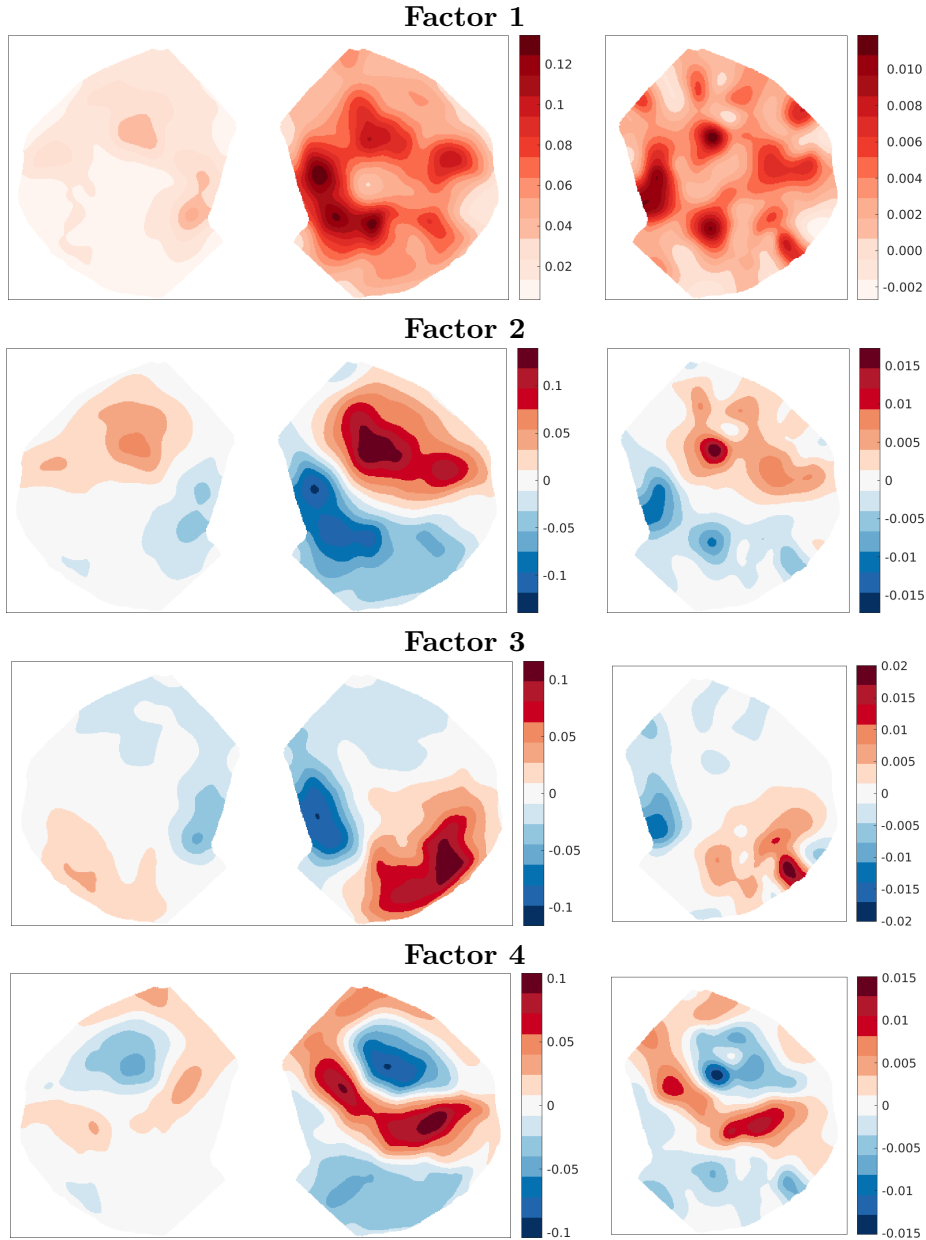

 \centering
 \captionof{figure}{Top four singular vectors of the flatmap connectivity
   with $\rank = 500$.
   Left: scaled target vectors $\hat U \Sigma$. Right: source vectors $\hat V$.}
 \label{fig:flatmap}
  \foreach \i in {1,...,4}{
    {\small \bf Factor \i} \\
    \vskip 1mm
    \includegraphics[height=3.5cm]{flat_US_\i.png}
    \hskip 3mm
    \includegraphics[height=3.5cm]{flat_V_\i.png}
    \\
  }
\end{figure}

Finally, we test the method on another problem which is a flatmap
projection of the brain (see Section~\ref{sec:data} for details).
This projection more faithfully represents areas of the cortex which
are missing from the top view since they curl underneath that vantage point.
The flatmap is closer to the kind of transformation used by cartographers
to flatten the globe, whereas the top view is like a satellite image
taken far from the globe.

The problem size is now larger by roughly a factor of three relative to
the top view.
Here, $\ny = 126\,847$ and $\nx = 63\,435$.
The number of experiments is the same, $\ninj=126$,
whereas the regularization parameter is set to $\bar\lambda=3 \times 10^7$.
The smoothing parameter was set to give the same level of smoothness,
measured ``by eye,'' in the components
as in the top view experiment.
The tolerances $\tau$ were as in the top view case.

In this case,
the computing time of the full solver would be excessively large,
so we do not estimate the error by comparison to the full solution,
instead taking the solution with $\rank = 1000$
as the reference solution $W_\star = W_{1000}$.
The computing times and the errors are shown in Figure~\ref{fig:flatmap_times}.
Here, the benefits by using the GPU implementation for
solving~\eqref{Proj_Sylv} were more significant than for the top view case.
We obtained the rank 500 solution in approximately 1.5 hours,
which is significantly less than pure CPU implementation,
which took 6.4 hours.
Comparing Figs.~\ref{fig:top_view_times} and \ref{fig:flatmap_times},
the computation times for the flatmap problem with $\rank = 500$ and 1000
are roughly twice as large as for the top view problem.
On the other hand, for $\rank = 125$ and 250, the compute times are
about three times as long for flatmap versus top view.
The observed scaling in compute time appears slightly slower than
$\mathcal{O}(n)$ in these tests.
The growth rate of the computing time on the GPU is
better than that of the CPU version since the matrix multiplications,
which dominate the CPU cost for large ranks,
are calculated in nearly constant time,
mainly due to communication overhead, on the GPU.
The RMS error between rank 500 and 1000
is again less than $10^{-6}$, so we believe rank 500
is probably a very good approximation to the full solution.
Figure~\ref{fig:J} shows the costs versus the approximation rank.
Again, we see that $\rank=500$ is reasonable and the distance from $W_\star$
is smaller than $10\%$.

The four dominant singular vectors of the
flatmap solution are shown in Figure~\ref{fig:flatmap},
oriented as in Figure~\ref{fig:top_view},
with the anterior-posterior axis from top-bottom and
the medial-lateral axis from left-right.
The first two factors are directly comparable between the
two problem outputs, although we see more structure in the
flatmap components.
This could be due to employing a projection which more accurately
represents these 3-D data in 2-D,
or due to the choice of smoothing parameter $\bar\lambda$.
The third and fourth components, on the other hand, are comparable
to the fourth and third components in the top view problem,
respectively.
Again, these patterns are reasonable and expected.
The raw 3-D data that were fed into the top view and
flatmap projections were the same,
but the greedy algorithm is run using different projected datasets.
It is reassuring that we can interpret the first few factors and
directly compare them against those in the top view.

\subsection{Dropping the nonnegativity constraint does not strongly affect the solutions}
In order to apply linear methods,
we relaxed the nonnegativity constraint when formulating
problem~\eqref{eq:obj_fun} (the unconstrained problem),
as opposed to the original problem~\eqref{eq:obj_fun_nonneg}
(the problem with nonnegativity constraint).
We now show that the resulting solutions are not significantly different
between the two problems.
This justifies the major simplification that we have made.

In all of our experiments with the test problem (Section~\ref{sec:test1}),
the resulting matrices were nearly nonnegative.
The solution $W_{40}$
has 48 out of 40\,000 negative entries.
These negative entries were all greater than -0.0023
in the lower-left corner of the matrix (see Figure~\ref{fig:test}),
where the truth is approximately zero.

We were able to solve the top view problem with the
nonnegative constraint using L-BFGS-B by initializing
with $W_{\rm full}$
projected onto the nonnegative orthant.
Let $W_{\rm proj}$ be the matrix with entries
$(W_{\rm proj})_{ij} = \max( 0, (W_\mathrm{full})_{ij})$,
and let $W_{\rm nonneg}$ denote the solution to the constrained
problem obtained in this way.
Comparing the nonnegative versus unconstrained solutions,
we found that
$\mathcal{E}(W_{\rm full}, W_{\rm nonneg}) = $ 3.99e-04.
Projecting $W_{\rm full}$ onto the nonnegative orthant leads to
$\mathcal{E}(W_{\rm proj}, W_{\rm nonneg}) = $ 3.67e-04.
In either case the $\infty$-norm difference is 0.009.
These results show that the solution to the unconstrained problem
is close to the solution of the constrained problem,
and that the projection of the solution to the unconstrained
problem is also close to the constrained solution.
Algorithm~\ref{alg:greedylr} thus offers an efficient way to
approximate the solution to the more difficult nonnegative problem,
while retaining low rank.

\section{Discussion}

\label{sec:discussion}

We have studied a numerical method
specifically tailored for the important neuroscience problem of
connectome regression from mesoscopic tract tracing experiments.
This connectome inference problem was formulated as
the regression problem~\eqref{eq:obj_fun}.
The optimality conditions for this problem turn out to be
a linear matrix equation in the unknown connectivity $W$,
which we propose to solve with Algorithm~\ref{alg:greedylr}.
Our numerical results show that the low-rank greedy algorithm,
as proposed by \cite{KreS15},
is a viable choice for acquiring low-rank factors of $W$
with a computation cost that was significantly
smaller compared to other approaches~\cite{Harris16,BenB13,KreT10}.
This allows us to infer the flatmap matrix, with approximately
$140\times$ more entries than previously obtained for the visual system,
while taking significantly less time: computing the flatmap solution
took hours versus days for the smaller low-rank visual network \cite{Harris16}.
The first few singular vector components of these cortical connectivities
are interpretable and reasonable from a neuroanatomy standpoint, although a
full anatomical study of this inferred connectivity is outside
the scope of the current paper.

The main ingredients of Algorithm~\ref{alg:greedylr} are
solving the large, sparse linear systems of equations at each
ALS iteration and solving the dense but small projected version
of the original linear matrix equation for the Galerkin step.
We had to carefully choose the solvers for each of these phases of the algorithm.
The Galerkin step forms the principal bottleneck
due to the absence of direct numerical methods to handle dense
linear matrix equations of moderate size.
We employed a matrix-valued CG iteration to approximately
solve~\eqref{Proj_Sylv},
implementing it on the GPU for speed.
This lead to cubic complexity
in $\rank$ at this step.
One could argue that equipping this CG iteration
with a preconditioner could speed up its convergence,
but so far we were not successful in finding a preconditioner
that both reduced the number of CG steps and the computational time.
A future research direction could be to derive an adequate preconditioning strategy
for the problem structure in~\eqref{multitermSylv},
that would increase the efficiency of any Krylov method.

Matrix-valued
Krylov subspace methods~\cite{Dam08,KreT10,BenB13,PalK19}
offer an alternative class of possible algorithms
to solving the overall linear matrix equation~\eqref{multitermSylv}.
However, for rapid convergence of these methods
we typically need a preconditioner.
In our tests on \eqref{multitermSylv},
these approaches performed poorly,
because rank truncations
(e.g.\ via thin QR or SVD)
are required after major
subcalculations which occur at every iteration.
Computing these decompositions quickly became
expensive because of the sheer amount of necessary rank truncations in the Krylov method.
If a suitable preconditioner for our problem would be found,
it would make sense to give low-rank matrix-valued Krylov
methods another try.

The original regression problem proposed by \cite{Harris16}
\eqref{eq:obj_fun_nonneg}
demands that the solution $W$ be nonnegative.
So far, this constraint is not considered by the employed algorithm.
However, for the test problem and data we have tried,
the computed matrix turns out to be majority
nonnegative.
We find typically small negative entries that can be safely neglected
without sacrificing accuracy.
Although a mostly nonnegative solution is not generally expected
when solving the unconstrained problem~\eqref{eq:obj_fun},
it appears that such behavior is typical for nonnegative
data matrices $X$ and $Y$.

Working directly with nonnegative factors $U \geq 0$ and $V \geq 0$
was originally proposed by
\cite{Harris16},
where they applied a projected gradient method to find such an
approximation for connectome of mouse visual areas albeit very slowly.
Such a formulation is preferred,
since it leads to a nonnegative $W$,
and it allows interpreting the leading factors as
the most important neural pathways in the brain.
Modifying Algorithm~\ref{alg:greedylr}
to compute nonnegative low-rank factors
or enforcing that the low-rank approximation
$UV^\trans\approx W$ is nonnegative---a nonlinear constraint---is a much harder goal to achieve.
For instance, even if one generated nonnegative factor matrices $U$ and $V$,
e.g.\ by changing the ALS step to nonnegative ALS,
the orthogonalization and Galerkin update each
destroy this nonnegativity.
New methods of NMF
which incorporate regularizations similar to our Laplacian
terms \cite{cichocki2009, cai2011}
are an area of ongoing research,
and the optimization techniques developed there could accelerate the
nonnegative low-rank formulation of \eqref{eq:obj_fun_nonneg}.
These include other techniques developed with neuroscience in mind,
such as neuron segmentation and calcium deconvolution
\cite{pnevmatikakis2016}
as well as sequence identification \cite{mackevicius2018}.
The greedy method we have presented is an excellent
way to initialize the nonnegative version of the problem,
similar to how SVD is used to initialize NMF.
We hope to improve upon nonnegative low-rank methods in the future.

Model \eqref{eq:obj_fun_nonneg}
is certainly not the only approach to solving
the connectome inference problem.
The loss term $\| P_\Omega (WX - Y) \|_F^2$
is standard and arises from Gaussian noise assumptions
combined with missing data and is
standard loss in matrix completion problems with noisy observations, e.g.,~\cite{mazumder2010,candes2010a}.
The regularization term is a thin plate spline penalty \cite{wahba1990}.
This is one of many
possible choices for smoothing,
among them penalties such as
$\| \mathrm{grad}(W) \|^2$ or the
total variation semi-norm \cite{rudin1992, chambolle2016},
which favors piecewise-constant solutions.
While we recognize that there are many possible choices
for the regularizer,
the thin plate penalty is reasonable,
linear and thus convenient to work with.
Previous work \cite{Harris16}
has shown that it is useful
for the connectome problem.
Testing other forms of regularization is a worthy goal
but not straightforward to implement at scale.
This is outside the scope of the current paper.

Finally, the most exciting prospects for this class of algorithms is what can
be learned when we apply them to next-generation tract tracing datasets.
Such techniques can be used to resolve differences
between the rat \cite{bota2003} brain and mouse \cite{Oh14},
or uncover unknown topographies, see~\cite{reimann2019}
in these and other animals (like the marmoset~\cite{majka2016}).
The mesoscale is also naturally the same resolution as obtained by wide-field calcium imaging.
Spatial connectome modeling could elucidate the largely
mysterious interactions different sensory modalities, proprioception, and motor areas, hopefully
leading to better understanding of integrative functions.

\section{Data and code}
\label{sec:data}
We tested our algorithm on two datasets
(top view and flatmap)
generated from Allen Institute for Brain Science
Mouse Connectivity Atlas data
\url{http://connectivity.brain-map.org}.
These data were obtained with the Python SDK
{\tt allensdk} version 0.13.1 available from
\url{http://alleninstitute.github.io/AllenSDK/}.
Our data pulling and processing scripts are available from
\url{https://github.com/kharris/allen-voxel-network}.

We used the {\tt allensdk} to retrieve
10 \micron injection and projection density volumetric data
for 126 {\it wildtype} experiments in cortex.
These data were projected from 3-D to 2-D
using either the top view or flatmap paths and saved as 2-D arrays.
Next, the projected coordinates were split into left and right
hemispheres.
Since {\it wildtype}
injections were always delivered into the right hemisphere,
this becomes our source space $S$ whereas the union of left and right
are the target space $T$.
We constructed 2-D 5-point Laplacian
matrices on these grids with ``free'' Neumann boundary conditions on the cortical edge.
Finally, the 2-D projected data were
downsampled 4 times along each dimension to obtain 40 \micron resolution.
The injection and projection data were then stacked into the matrices
$X$ and $Y$, respectively.
The mask $\Omega$ was set via $\Omega_{ij} = 1_{\{ X_{ij} \leq 0.4 \}}$.

\matlab~code which implements our greedy low-rank algorithm
\eqref{alg:greedylr} is included in the repository:
\url{https://gitlab.mpi-magdeburg.mpg.de/kuerschner/lowrank\_connectome}.
We also include the problem inputs $X, Y, L_x, L_y, \Omega$
for our three example problems
(test, top view, and flatmap) as \matlab~files.
Note that $\Omega$ is stored as $1-\Omega$ in these files,
as this matches the convention of \cite{Harris16}.


\section*{Acknowledgements}
We would like to thank Lydia Ng, Nathan Gouwens, Stefan Mihalas,
Nile Graddis and others at the Allen Institute
for the top view and flatmap paths and general help accessing the data.
Thank you to Braden Brinkman for discussions of the continuous problem,
to Stefan Mihalas and Eric Shea-Brown for general discussions.
This work was primarily generated while PK was affiliated with the
Max Planck Institute for Dynamics of Complex Technical Systems.

KDH was supported by the
Big Data for Genomics and Neuroscience NIH training grant
and a Washington Research Foundation Postdoctoral Fellowship.
SD is thankful to the Engineering and Physical Sciences Research Council (UK)
for supporting his postdoctoral position at the University of Bath through Fellowship EP/M019004/1,
and the kind hospitality of the Erwin Schr\"odinger International Institute for Mathematics and Physics (ESI),
where this manuscript was finalised during the Thematic Programme
{\it Numerical Analysis of Complex PDE Models in the Sciences}.


\begin{thebibliography}{10}

\bibitem{altas1998}
{\sc I.~Altas, J.~Dym, M.~Gupta, and R.~Manohar}, {\em Multigrid {{Solution}}
  of {{Automatically Generated High}}-{{Order Discretizations}} for the
  {{Biharmonic Equation}}}, SIAM Journal on Scientific Computing, 19 (1998),
  pp.~1575--1585, \url{https://doi.org/10.1137/S1464827596296970}.

\bibitem{Ben04}
{\sc P.~Benner}, {\em Solving large-scale control problems}, {IEEE} Control
  Syst. Mag., 14 (2004), pp.~44--59.

\bibitem{BenB13}
{\sc P.~Benner and T.~Breiten}, {\em {Low rank methods for a class of
  generalized {L}yapunov equations and related issues}}, {Numer. Math.}, 124
  (2013), pp.~441--470, \url{https://doi.org/10.1007/s00211-013-0521-0}.

\bibitem{BenLT09}
{\sc P.~Benner, R.-C. Li, and N.~Truhar}, {\em On the {ADI} method for
  {S}ylvester equations}, J. Comput. Appl. Math., 233 (2009), pp.~1035--1045.

\bibitem{BenS13}
{\sc P.~Benner and J.~Saak}, {\em Numerical solution of large and sparse
  continuous time algebraic matrix {R}iccati and {L}yapunov equations: a state
  of the art survey}, GAMM Mitteilungen, 36 (2013), pp.~32--52,
  \url{https://doi.org/10.1002/gamm.201310003}.

\bibitem{bota2003}
{\sc M.~Bota, H.-W. Dong, and L.~W. Swanson}, {\em {From Gene Networks to Brain
  Networks}}, Nature Neuroscience, 6 (2003), pp.~795--799,
  \url{https://doi.org/10.1038/nn1096}.

\bibitem{buckner2019}
{\sc R.~L. Buckner and D.~S. Margulies}, {\em Macroscale cortical organization
  and a default-like apex transmodal network in the marmoset monkey}, Nature
  Communications, 10 (2019), \url{https://doi.org/10.1038/s41467-019-09812-8}.

\bibitem{byrd1995}
{\sc R.~Byrd, P.~Lu, J.~Nocedal, and C.~Zhu}, {\em A {{Limited Memory
  Algorithm}} for {{Bound Constrained Optimization}}}, SIAM Journal on
  Scientific Computing, 16 (1995), pp.~1190--1208,
  \url{https://doi.org/10.1137/0916069}.

\bibitem{cai2011}
{\sc D.~Cai, X.~He, J.~Han, and T.~S. Huang}, {\em {Graph {{Regularized
  Nonnegative Matrix Factorization}} for {{Data Representation}}}}, IEEE
  Transactions on Pattern Analysis and Machine Intelligence, 33 (2011),
  pp.~1548--1560, \url{https://doi.org/10.1109/TPAMI.2010.231}.

\bibitem{candes2010a}
{\sc E.~J. Candes and Y.~Plan}, {\em Matrix {{Completion With Noise}}},
  Proceedings of the IEEE, 98 (2010), pp.~925--936,
  \url{https://doi.org/10.1109/JPROC.2009.2035722}.

\bibitem{chambolle2016}
{\sc A.~Chambolle and T.~Pock}, {\em An introduction to continuous optimization
  for imaging}, Acta Numerica, 25 (2016), pp.~161--319,
  \url{https://doi.org/10.1017/S096249291600009X}.

\bibitem{cichocki2009}
{\sc A.~Cichocki, R.~Zdunek, A.~H. Phan, and S.-i. Amari}, {\em {Nonnegative
  {{Matrix}} and {{Tensor Factorizations}}: {{Applications}} to {{Exploratory
  Multi}}-Way {{Data Analysis}} and {{Blind Source Separation}}}}, {John Wiley
  \& Sons}, July 2009.

\bibitem{Dam08}
{\sc T.~Damm}, {\em {Direct methods and {ADI}-preconditioned {K}rylov subspace
  methods for generalized {L}yapunov equations}}, {Numer. Lin. Alg. Appl.}, 15
  (2008), pp.~853--871.

\bibitem{JarKMR16}
{\sc {E. Ringh, G. Mele, J. Karlsson, E. Jarlebring}}, {\em {Sylvester-based
  preconditioning for the waveguide eigenvalue problem}}, {Linear Algebra
  Appl.}, 542 (2018), pp.~441--463.

\bibitem{GolV13}
{\sc G.~H. Golub and C.~F. {Van~Loan}}, {\em {Matrix Computations}}, Johns
  Hopkins University Press, Baltimore, fourth~ed., 2013.

\bibitem{Gra04}
{\sc L.~Grasedyck}, {\em {Existence of a low rank or {$H$}-matrix approximant
  to the solution of a {S}ylvester equation}}, {Numer. Lin. Alg. Appl.}, 11
  (2004), pp.~371--389.

\bibitem{grillner2016}
{\sc S.~Grillner, N.~Ip, C.~Koch, W.~Koroshetz, H.~Okano, M.~Polachek, M.-m.
  Poo, and T.~J. Sejnowski}, {\em {Worldwide Initiatives to Advance Brain
  Research}}, Nature Neuroscience,  (2016),
  \url{https://doi.org/10.1038/nn.4371}.

\bibitem{gamanut2018}
{\sc R.~{G\u{a}m\u{a}nu{\c t}}, H.~Kennedy, Z.~Toroczkai, M.~{Ercsey-Ravasz},
  D.~C. {Van Essen}, K.~Knoblauch, and A.~Burkhalter}, {\em {The {{Mouse
  Cortical Connectome}}, {{Characterized}} by an {{Ultra}}-{{Dense Cortical
  Graph}}, {{Maintains Specificity}} by {{Distinct Connectivity Profiles}}}},
  Neuron, 97 (2018), pp.~698--715.e10,
  \url{https://doi.org/10.1016/j.neuron.2017.12.037}.

\bibitem{HardtW14}
{\sc M.~Hardt and M.~Wootters}, {\em Fast matrix completion without the
  condition number}, in Proceedings of The 27th Conference on Learning Theory,
  {COLT} 2014, Barcelona, Spain, June 13-15, 2014, 2014, pp.~638--678.

\bibitem{Harris16}
{\sc K.~D. Harris, S.~Mihalas, and E.~Shea-Brown}, {\em {High Resolution Neural
  Connectivity from Incomplete Tracing Data Using Nonnegative Spline
  Regression}}, in {Neural {{Information Processing Systems}}}, 2016.

\bibitem{harshman70}
{\sc R.~Harshman}, {\em Foundations of the {PARAFAC} procedure: Models and
  conditions for an ``explanatory'' multi-modal factor analysis}, UCLA Working
  Papers in Phonetics, 16 (1970).

\bibitem{JarMPR17}
{\sc E.~{Jarlebring}, G.~{Mele}, D.~{Palitta}, and E.~{Ringh}}, {\em {Krylov
  methods for low-rank commuting generalized Sylvester equations}}, {Numer.
  Lin. Alg. Appl.},  (2018), \url{https://doi.org/10.1002/nla.2176}.

\bibitem{jenett2012}
{\sc A.~Jenett, G.~M. Rubin, T.-T.~B. Ngo, D.~Shepherd, C.~Murphy, H.~Dionne,
  B.~D. Pfeiffer, A.~Cavallaro, D.~Hall, J.~Jeter, N.~Iyer, D.~Fetter, J.~H.
  Hausenfluck, H.~Peng, E.~T. Trautman, R.~R. Svirskas, E.~W. Myers, Z.~R.
  Iwinski, Y.~Aso, G.~M. DePasquale, A.~Enos, P.~Hulamm, S.~C.~B. Lam, H.-H.
  Li, T.~R. Laverty, F.~Long, L.~Qu, S.~D. Murphy, K.~Rokicki, T.~Safford,
  K.~Shaw, J.~H. Simpson, A.~Sowell, S.~Tae, Y.~Yu, and C.~T. Zugates}, {\em {A
  {{GAL4}}-{{Driver Line Resource}} for {{Drosophila Neurobiology}}}}, Cell
  Reports, 2 (2012), pp.~991--1001,
  \url{https://doi.org/10.1016/j.celrep.2012.09.011}.

\bibitem{kasthuri2015}
{\sc N.~Kasthuri, K.~J. Hayworth, D.~R. Berger, R.~L. Schalek, J.~A. Conchello,
  S.~{Knowles-Barley}, D.~Lee, A.~{V{\'a}zquez-Reina}, V.~Kaynig, T.~R. Jones,
  M.~Roberts, J.~L. Morgan, J.~C. Tapia, H.~S. Seung, W.~G. Roncal, J.~T.
  Vogelstein, R.~Burns, D.~L. Sussman, C.~E. Priebe, H.~Pfister, and J.~W.
  Lichtman}, {\em Saturated {{Reconstruction}} of a {{Volume}} of
  {{Neocortex}}}, Cell, 162 (2015), pp.~648--661,
  \url{https://doi.org/10.1016/j.cell.2015.06.054}.

\bibitem{kennedy2016}
{\em Micro-, {{Meso}}- and {{Macro}}-{{Connectomics}} of the {{Brain}}},
  Research and {{Perspectives}} in {{Neurosciences}}, {Springer International
  Publishing}, 2016.

\bibitem{Knox18}
{\sc J.~E. Knox, K.~D. Harris, N.~Graddis, J.~D. Whitesell, H.~Zeng, J.~A.
  Harris, E.~Shea-Brown, and S.~Mihalas}, {\em {High Resolution Data-Driven
  Model of the Mouse Connectome}}, bioRxiv,  (2018), p.~293019,
  \url{https://doi.org/10.1101/293019}.

\bibitem{KreS15}
{\sc D.~Kressner and P.~Sirkovi{\'c}}, {\em {Truncated low-rank methods for
  solving general linear matrix equations}}, {Numer. Lin. Alg. Appl.}, 22
  (2015), pp.~564--583, \url{https://doi.org/10.1002/nla.1973}.

\bibitem{KreT10}
{\sc D.~Kressner and C.~Tobler}, {\em {Krylov Subspace Methods for Linear
  Systems with Tensor Product Structure}}, {{SIAM} J. Matrix Anal. Appl.}, 31
  (2010), pp.~1688--1714.

\bibitem{lein2007}
{\sc E.~S. Lein, M.~J. Hawrylycz, N.~Ao, M.~Ayres, A.~Bensinger, A.~Bernard,
  A.~F. Boe, M.~S. Boguski, K.~S. Brockway, E.~J. Byrnes, L.~Chen, L.~Chen,
  T.-M. Chen, M.~C. Chin, J.~Chong, B.~E. Crook, A.~Czaplinska, C.~N. Dang,
  S.~Datta, N.~R. Dee, A.~L. Desaki, T.~Desta, E.~Diep, T.~A. Dolbeare, M.~J.
  Donelan, H.-W. Dong, J.~G. Dougherty, B.~J. Duncan, A.~J. Ebbert, G.~Eichele,
  L.~K. Estin, C.~Faber, B.~A. Facer, R.~Fields, S.~R. Fischer, T.~P. Fliss,
  C.~Frensley, S.~N. Gates, K.~J. Glattfelder, K.~R. Halverson, M.~R. Hart,
  J.~G. Hohmann, M.~P. Howell, D.~P. Jeung, R.~A. Johnson, P.~T. Karr,
  R.~Kawal, J.~M. Kidney, R.~H. Knapik, C.~L. Kuan, J.~H. Lake, A.~R. Laramee,
  K.~D. Larsen, C.~Lau, T.~A. Lemon, A.~J. Liang, Y.~Liu, L.~T. Luong,
  J.~Michaels, J.~J. Morgan, R.~J. Morgan, M.~T. Mortrud, N.~F. Mosqueda, L.~L.
  Ng, R.~Ng, G.~J. Orta, C.~C. Overly, T.~H. Pak, S.~E. Parry, S.~D. Pathak,
  O.~C. Pearson, R.~B. Puchalski, Z.~L. Riley, H.~R. Rockett, S.~A. Rowland,
  J.~J. Royall, M.~J. Ruiz, N.~R. Sarno, K.~Schaffnit, N.~V. Shapovalova,
  T.~Sivisay, C.~R. Slaughterbeck, S.~C. Smith, K.~A. Smith, B.~I. Smith, A.~J.
  Sodt, N.~N. Stewart, K.-R. Stumpf, S.~M. Sunkin, M.~Sutram, A.~Tam, C.~D.
  Teemer, C.~Thaller, C.~L. Thompson, L.~R. Varnam, A.~Visel, R.~M. Whitlock,
  P.~E. Wohnoutka, C.~K. Wolkey, V.~Y. Wong, M.~Wood, M.~B. Yaylaoglu, R.~C.
  Young, B.~L. Youngstrom, X.~F. Yuan, B.~Zhang, T.~A. Zwingman, and A.~R.
  Jones}, {\em {Genome-Wide Atlas of Gene Expression in the Adult Mouse
  Brain}}, Nature, 445 (2007), pp.~168--176,
  \url{https://doi.org/10.1038/nature05453}.

\bibitem{mackevicius2018}
{\sc E.~L. Mackevicius, A.~H. Bahle, A.~H. Williams, S.~Gu, N.~I. Denissenko,
  M.~S. Goldman, and M.~S. Fee}, {\em {Unsupervised Discovery of Temporal
  Sequences in High-Dimensional Datasets, with Applications to Neuroscience}},
  bioRxiv,  (2018), p.~273128, \url{https://doi.org/10.1101/273128}.

\bibitem{majka2016}
{\sc P.~Majka, {Chaplin, Tristan A.}, {Yu, Hsin-Hao}, {Tolpygo, Alexander},
  {Mitra, Partha P.}, {W{\'o}jcik, Daniel K.}, and {Rosa, Marcello G.P.}}, {\em
  {Towards a Comprehensive Atlas of Cortical Connections in a Primate Brain:
  {{Mapping}} Tracer Injection Studies of the Common Marmoset into a Reference
  Digital Template}}, Journal of Comparative Neurology, 524 (2016),
  pp.~2161--2181, \url{https://doi.org/10.1002/cne.24023}.

\bibitem{mazumder2010}
{\sc R.~Mazumder, T.~Hastie, and R.~Tibshirani}, {\em Spectral {{Regularization
  Algorithms}} for {{Learning Large Incomplete Matrices}}}, J. Mach. Learn.
  Res., 11 (2010), pp.~2287--2322.

\bibitem{mitra2014}
{\sc P.~P. Mitra}, {\em {The {{Circuit Architecture}} of {{Whole Brains}} at
  the {{Mesoscopic Scale}}}}, Neuron, 83 (2014), pp.~1273--1283,
  \url{https://doi.org/10.1016/j.neuron.2014.08.055}.

\bibitem{nouy-greedy-2010}
{\sc A.~Nouy}, {\em Proper generalized decompositions and separated
  representations for the numerical solution of high dimensional stochastic
  problems}, Archives of Computational Methods in Engineering, 17 (2010),
  pp.~403--434, \url{https://doi.org/10.1007/s11831-010-9054-1}.

\bibitem{Oh14}
{\sc S.~W. Oh, J.~A. Harris, L.~Ng, B.~Winslow, N.~Cain, S.~Mihalas, Q.~Wang,
  C.~Lau, L.~Kuan, A.~M. Henry, M.~T. Mortrud, B.~Ouellette, T.~N. Nguyen,
  S.~A. Sorensen, C.~R. Slaughterbeck, W.~Wakeman, Y.~Li, D.~Feng, A.~Ho,
  E.~Nicholas, K.~E. Hirokawa, P.~Bohn, K.~M. Joines, H.~Peng, M.~J. Hawrylycz,
  J.~W. Phillips, J.~G. Hohmann, P.~Wohnoutka, C.~R. Gerfen, C.~Koch,
  A.~Bernard, C.~Dang, A.~R. Jones, and H.~Zeng}, {\em {A Mesoscale Connectome
  of the Mouse Brain}}, Nature, 508 (2014), pp.~207--214,
  \url{https://doi.org/10.1038/nature13186}.

\bibitem{OrtR00}
{\sc J.~M. Ortega and W.~C. Rheinboldt}, {\em {Iterative solution of nonlinear
  equations in several variables}}, SIAM, 2000.

\bibitem{PalK19}
{\sc D.~Palitta and P.~K{\"u}rschner}, {\em On the convergence of Krylov
  methods with low-rank truncations}, e-print 1909.01226, arXiv, 2019,
  \url{https://arxiv.org/abs/1909.01226}.
\newblock math.NA.

\bibitem{pati-omp-1993}
{\sc Y.~C. {Pati}, R.~{Rezaiifar}, and P.~S. {Krishnaprasad}}, {\em Orthogonal
  matching pursuit: recursive function approximation with applications to
  wavelet decomposition}, in Proceedings of 27th Asilomar Conference on
  Signals, Systems and Computers, vol.~1, 1993, pp.~40--44,
  \url{https://doi.org/10.1109/ACSSC.1993.342465}.

\bibitem{pnevmatikakis2016}
{\sc E.~A. Pnevmatikakis, D.~Soudry, Y.~Gao, T.~A. Machado, J.~Merel, D.~Pfau,
  T.~Reardon, Y.~Mu, C.~Lacefield, W.~Yang, M.~Ahrens, R.~Bruno, T.~M. Jessell,
  D.~S. Peterka, R.~Yuste, and L.~Paninski}, {\em {Simultaneous {{Denoising}},
  {{Deconvolution}}, and {{Demixing}} of {{Calcium Imaging Data}}}}, Neuron, 89
  (2016), pp.~285--299, \url{https://doi.org/10.1016/j.neuron.2015.11.037}.

\bibitem{PowSS17}
{\sc C.~E. Powell, D.~Silvester, and V.~Simoncini}, {\em {An Efficient Reduced
  Basis Solver for Stochastic Galerkin Matrix Equations}}, {{SIAM} J. Sci.
  Comput.}, 39 (2017), pp.~A141--A163,
  \url{https://doi.org/10.1137/15M1032399}.

\bibitem{reimann2019}
{\sc M.~W. Reimann, M.~Gevaert, Y.~Shi, H.~Lu, H.~Markram, and E.~Muller}, {\em
  A null model of the mouse whole-neocortex micro-connectome}, Nature
  Communications, 10 (2019), pp.~1--16,
  \url{https://doi.org/10.1038/s41467-019-11630-x}.

\bibitem{rudin1992}
{\sc L.~I. Rudin, S.~Osher, and E.~Fatemi}, {\em Nonlinear total variation
  based noise removal algorithms}, Physica D: Nonlinear Phenomena, 60 (1992),
  pp.~259--268, \url{https://doi.org/10.1016/0167-2789(92)90242-F}.

\bibitem{ShaSS15}
{\sc S.~D. Shank, V.~Simoncini, and D.~B. Szyld}, {\em {Efficient low-rank
  solution of generalized {L}yapunov equations}}, {Numer. Math.}, 134 (2015),
  pp.~327--342.

\bibitem{Sim16}
{\sc V.~Simoncini}, {\em {Computational methods for linear matrix equations}},
  {{SIAM} Rev.}, 38 (2016), pp.~377--441.

\bibitem{SorVBDeL13}
{\sc L.~Sorber, M.~Van~Barel, and L.~De~Lathauwer}, {\em Optimization-based
  algorithms for tensor decompositions: Canonical polyadic decomposition,
  decomposition in rank-$({L}_r,{L}_r,1)$ terms, and a new generalization},
  SIAM Journal on Optimization, 23 (2013), pp.~695--720,
  \url{https://doi.org/10.1137/120868323}.

\bibitem{sporns2010}
{\sc O.~Sporns}, {\em {Networks of the {{Brain}}}}, {The MIT Press}, 1st~ed.,
  2010.

\bibitem{vanessen2013}
{\sc D.~C. {Van Essen}}, {\em {Cartography and {{Connectomes}}}}, Neuron, 80
  (2013), pp.~775--790, \url{https://doi.org/10.1016/j.neuron.2013.10.027}.

\bibitem{wahba1990}
{\sc G.~Wahba}, {\em {Spline {{Models}} for {{Observational Data}}}}, {SIAM},
  Sept. 1990.

\bibitem{ypma2016}
{\sc R.~J.~F. Ypma and E.~T. Bullmore}, {\em {Statistical {{Analysis}} of
  {{Tract}}-{{Tracing Experiments Demonstrates}} a {{Dense}}, {{Complex
  Cortical Network}} in the {{Mouse}}}}, PLOS Comput Biol, 12 (2016),
  p.~e1005104, \url{https://doi.org/10.1371/journal.pcbi.1005104}.

\end{thebibliography}

\end{document}